\documentstyle{article}
\newtheorem{theorem}{Theorem}[section]
\newtheorem{lemma}[theorem]{Lemma}
\newtheorem{proposition}[theorem]{Proposition}
\newtheorem{corollary}[theorem]{Corollary}
\newtheorem{definition}[theorem]{Definition}

\begin{document}

\title{{\large {\bf VON NEUMANN SPECTRA NEAR THE SPECTRAL GAP}}}
\author{{Alan L. Carey} \\
{\small {\em Department of Pure Mathematics, University of Adelaide}} \\
{\small {\em Adelaide, South Australia, Australia.}}\\\\
{Thierry Coulhon}\\
{\small {\em D\'epartement de Math\'ematiques, Universit\'e de Cergy-Pontoise
}}\\
{\small {\em Cergy-Pontoise, France.}}\\\\
{Varghese Mathai}\\
{\small {\em Department of Pure Mathematics, University of Adelaide}} \\
{\small {\em Adelaide, South Australia, Australia.}}\\\\
{John Phillips}\\
{\small {\em Department of Mathematics, University of Victoria}}\\
{\small {\em Victoria, BC, Canada.}}}
\date{}
\maketitle

\noindent{\bf Keywords}. Heat kernels, bottom of the spectrum,
Novikov-Shubin invariants, von Neumann determinants, torsion.

\noindent{\bf AMS 1991 Subject Classification}. 58G11, 58G18 and 58G25.

\noindent{\bf Abstract}. In this paper we study some new von Neumann spectral
invariants associated to the Laplacian acting on $L^2$ differential forms on
the universal cover of a closed manifold. These invariants coincide with the
Novikov-Shubin invariants whenever there is no spectral gap in the spectrum
of the Laplacian, and are homotopy invariants in this case. In the presence
of a spectral gap, they differ in character and value from the
Novikov-Shubin invariants. Under a positivity assumption on these invariants,
we prove that certain $L^2$ theta and $L^2$ zeta functions
defined by metric dependent combinatorial Laplacians acting on $L^2$ cochains
associated with a triangulation of the manifold, converge uniformly to
their analytic
counterparts, as the mesh of the triangulation goes to zero.

\vspace{.2in}

\noindent{\Large {\bf Introduction}}

Let $(M,g)$ be a compact connected Riemannian manifold and $\widetilde{M}$ be
its universal cover, which is a Riemannian manifold with the induced metric.
Although we only consider the universal cover in this paper, our results extend
easily to any Galois covering of $M$. Let $\Delta_{j} = d\delta + \delta d$
denote the Laplacian acting on $L^{2}$ differential $j$-forms on
$\widetilde{M}$, where $d$ denotes the de Rham differential and $\delta$ its
adjoint. Let $\lambda_{0,j}$ denote the bottom of the spectrum of $\Delta_{j}$
on the orthogonal complement of its kernel.

Let $k_{j}(t,x,y)$ denote the integral kernel of the heat operator $
e^{-t\Delta_ {j}}$. By parabolic regularity theory, $k_{j}(t,x,y)$ is smooth.
Since $\Delta_{j}$ commutes with the action of the fundamental group
$\pi_{1}(M)$, we see that \[ k_{j}(t,\gamma.{x},\gamma.{y})=k_{j}(t,x,y) \] for
all $\gamma\in\pi_{1}(M)$ (here we identify the cotangent spaces of $x$ and
$\gamma.x$ via $\gamma$). $\;$Recall that the von Neumann trace of $
e^{-t\Delta_j}$ is given by \[
\tau(e^{-t\Delta_{j}})=\int_{M}tr(k_{j}(t,x,x))dx. \] Atiyah \cite{A} defined
the $L^{2}$ Betti numbers of $M$ as follows, \[ b^{j}_{(2)}(\widetilde M) =
\tau(P_j), \] where $P_j$ denotes the orthogonal projection onto the kernel of
$\Delta_j$. Then, since $k_{j}(t,x,y)$ converges as $t\rightarrow \infty$
uniformly over compact subsets to $k_{P_j}(x,y)$, where $k_{P_j}(x,y)$ denotes
the integral kernel of the operator $P_j$ (cf.  \cite{R}), we see that \[
b^{j}_{(2)}(\widetilde M)=\lim_{t\rightarrow\infty}\tau(e^ {-t\Delta_{j}}). \]
Define the von Neumann algebra theta functions \[
\theta_{j}(t)=\tau(e^{-t\Delta_{j}}) - b^{j}_{(2)}(\widetilde M) \] Explicit
calculations tend to show that the large time asymptotics of $ \theta_{j}(t)$
are of the form $e^{-\lambda_{0,j}t}t^{-\beta}$ for some $ \beta>0$. This
motivates the following definitions \[ \beta_{j}(M,g)=\sup\{\beta\in
I\!\!R:e^{\lambda_{0,j}t}\theta_{j}(t)\;\; \mbox{ is } \;\;O(t^{-\beta})\;\mbox{
as }\;t\rightarrow\infty\}\in[0,\infty ] \] \[
\overline\beta_{j}(M,g)=\inf\{\beta\in I\!\!R: t^{-\beta}\;\;\mbox{ is }
\;\;O(e^{\lambda_{0,j}t}\theta_{j}(t)))\;\mbox{ as }\;t\rightarrow\infty\}
\in[0,\infty]. \] Whenever $\lambda_{0,j}=0$ (note that this condition is
independent of the choice of metric cf. \cite{GS}), $\beta_{j}(M,g)$ and
$\overline \beta_{j}(M,g)$ coincide with the Novikov-Shubin invariants (cf.
\cite{ES}, \cite{GS}, \cite{BMW}, \cite{L}) of $M$ and therefore they are
independent of the choice of metric. However when $\lambda_{0,j}>0$ general,
they differ from the Novikov-Shubin invariants,  and are probably dependent on
the choice of Riemannian metric, as an example by Donald Cartwright \cite{Ca}
appears to indicate. But $ \beta_{j}(M,g)$ and $\overline\beta_{j}(M,g)$ are by
fiat clearly von Neumann spectral invariants, which tend to be finite even when
Novikov-Shubin invariants aren't, as we demonstrate by examples in section 5.
They are also invariant under scaling of the metric by a constant. The Hodge
star operator on the universal cover of $M$ is a von Neumann algebra isometry
intertwining the Laplacians $\Delta_j$ and $\Delta_{n-j}$, where $n$ denotes the
dimension of $M$. Hence we see that $\beta_{j}(M,g) = \beta_{n-j}(M,g)$ and $
\overline\beta_{j}(M,g)=\overline\beta_{n-j}(M,g)$ for all $0\leq j\leq n$.
Using the local Harnack inequality \cite{VSC}, one observes that the large time
asymptotics of $\theta_{0}(t)$ coincide with the large time asymptotics of
$||e^{-t\Delta_{0}}||_{1\rightarrow\infty}$. It is unclear to us if this is also
true on $j$-forms, where $j>0$.

We next give a summary of our results. Let $K$ denotes a triangulation of $M$.
In section 1, we define a metric dependent combinatorial Laplacian
$\Delta_j^{\widetilde K, \widetilde W}$ acting on $L^2$ cochains on the
universal cover $\widetilde K$, besides other background material. Section 2
contains a detailed analysis of the relationship between the spectral density
functions for the combinatorial Laplacian and the analytic Laplacian. These
estimates are the analogues on non-compact covering spaces of some results of
Dodziuk and Patodi \cite{DP} that also enable us to show that the bottom of the
spectrum of the combinatorial Laplacian, $\lambda_{0,j}^{\widetilde K,
\widetilde W}$ converges to the bottom of the spectrum of the analytic Laplacian
$\lambda_{0,j}$, as the mesh of the triangulation goes to zero. In section 3 we
prove some estimates for the combinatorial and the analytic $L^2$ theta
functions, which are used later on in the paper. In section 4, we prove that
$e^{\lambda_{0,j}^{\widetilde K, \widetilde W}t} \theta_j^{\widetilde K,
\widetilde W}(t)$ converges uniformly in $t\in [t_0, \infty)$ to
$e^{\lambda_{0,j}t}\theta_j(t)$, for any $t_0>0$, as the mesh $\eta$ of the
triangulation $K$ goes to zero. In section 5, we compute the values of the
invariants $\beta_j(M,g), \bar\beta_j(M,g)$ for closed odd dimensional
hyperbolic manifolds. In section 6, under the assumption that $\beta_j(M,g)>0$,
we define the von Neumann determinant of the operator $\Delta_j-\lambda_{0,j}$
as well as the analytic $\beta$ torsion, which is defined analogously to
Ray-Singer torsion \cite{RS}, as an alternating product of the von Neumann
determinants of $\Delta_j-\lambda_{0,j}$. We also define its combinatorial
counterpart, which we call combinatorial $\beta$-torsion. In section 7, we prove
convergence theorems for $L^2$ theta and $L^2$ zeta functions respectively, as
the mesh of the triangulation goes to zero, extending to the non-compact
covering space classical theorems \cite {DP}. The latter arise naturally in the
study of $L^2$ torsion invariants, see \cite{M1,CM} and \cite{L}. We anticipate
that the technical results of this paper have applications in a number of
directions and will provide the groundwork for the investigation of other
spectral properties of the Laplacian on non-compact manifolds. The results of
section 7 also give evidence that the combinatorial $\beta$-torsion converges to
the analytic $\beta$-torsion, as the mesh of the triangulation goes to zero. We
plan to study this and other questions raised in this paper elsewhere.

We conclude the introduction with the conjecture that $\beta_j(M,g)$ is always
positive. Earlier, Sunada \cite{Su2} had conjectured that this was the case for
functions, that is, $\beta_0(M,g)>0$. This has been proved by Terry Lyons and we
present his proof in the appendix to our paper. We also conjecture that
$\beta_j(K,g)$ converges to $\beta_j(M,g)$ as the mesh of the triangulation goes
to zero, where $\beta_j(K,g)$ denotes the combinatorial counterpart of
$\beta_j(M,g)$. The results of section 4 suggest that this conjecture is true,
and one can find a discussion of it over there. \bigskip

\section{ Preliminaries.}

In this section we establish the notation of the paper. We recall the definition
of the combinatorial and analytic $L^2$ theta functions and also some of their
basic properties  following the notation established in \cite {CM} and \cite{Do}
(with some minor changes). Thus $(M,g)$ is a compact manifold without boundary
with Riemannian metric $g$ and dimension $n$. The fundamental group of $M$ we
denote by $\Gamma$ and we lift $g$ to a $\Gamma$ invariant Riemannian metric
$\tilde g$ on the universal cover $\widetilde M$ . Now introduce the space
$\Omega_{(2)}(\widetilde M)$ of $L^2$ differential forms with respect to the
volume defined by this metric. Denote by $d$ the exterior derivative (which is a
$\Gamma$-invariant lift of the exterior derivative on forms on $M$) on
$\Omega_{(2)}(\widetilde M)$. Strictly speaking we must restrict the action of
$d$ to sufficiently smooth forms. When this becomes crucial we shall introduce
the appropriate notation. Let $K $ be a triangulation of $M$ and $\widetilde K$
be the induced $\Gamma$ -invariant triangulation of $\widetilde M$. Let $C(K)$
denote the space of cochains on $K$ and $C_{(2)}(\widetilde K)$ denote the space
of $\ell^2$ cochains on $\widetilde K$.

Following \cite{Do}, we use the de Rham and Whitney maps to relate the
combinatorial and de Rham complexes. For $\widetilde M$, let $\widetilde A$
denote the de Rham map from the de Rham complex to the combinatorial complex and
$\widetilde W$ denote the Whitney map which goes in the reverse direction and
satisfies: \[ \widetilde A\widetilde W = 1 \] (we will review the definitions of
these maps in section 2). We may use the Whitney map to define an equivalent
inner product on $C_{(2)}(\widetilde K)$ by setting \[ <c_\sigma,
c_{\sigma^\prime}>_{\widetilde W} = <\widetilde Wc_\sigma, \widetilde
Wc_{\sigma^\prime}> \] where $\sigma,\sigma^\prime$ are elements of $\widetilde
K$ and $c_\sigma$ denotes the characteristic cochain of $\sigma$. The
fundamental group $\Gamma $ acts on each of the complexes
$\Omega_{(2)}(\widetilde M)$ and $ C_{(2)}(\widetilde K)$ by unitary operators.
There are isomorphisms of $ \Gamma$-modules: \[ \Omega_{(2)}(\widetilde
M)\cong\Omega_{(2)}(M)\otimes \ell^2(\Gamma),\ \ \ C_{(2)}(\widetilde K)\cong
C(K)\otimes \ell^2(\Gamma), \] where $\Omega_{(2)}(M)$ denotes $L^2$ forms on
$M$ and $\Gamma$ acts on $ \ell ^2(\Gamma)$ by the left regular representation
$\lambda$. The von Neumann algebra generated by $\{\lambda(\gamma)\ \vert\
\gamma\in \Gamma\}$ is a finite von Neumann algebra ${\cal U}$ with normalised
trace denoted $ \tau$. We may regard the trace $\tau$ as being defined on the
commutant as well (as in \cite{Di}) by letting $\chi_e$ denote the function in $
\ell^2(\Gamma)$ which is one at the identity of $\Gamma$ and zero elsewhere and
then $\tau(B)=<B\chi_e,\chi_e>_{\ell^2(\Gamma)}$, for any $B$ in ${\cal U }$ or
${\cal U}^{\prime}$. There is also a finite trace on the commutant $ {\cal
B}(C(K))\otimes {\cal U}^\prime$ of this ${\cal U}$ action on $
C_{(2)}(\widetilde K)$: it is ${\rm tr}\otimes\tau$ where ${\rm tr}$ denotes the
usual matrix trace on the finite dimensional space ${\cal B}(C(K))$. Note that
${\cal U}$ and ${\cal U}^\prime$ are anti-isomorphic \cite{Di}.

There is a semifinite trace also denoted tr$\otimes \tau$ on the commutant of
the ${\cal U}$ action on $\Omega_{(2)}(\widetilde M)$. This commutant is just
${\cal B}(\Omega_{(2)}(M))\otimes{\cal U}^\prime$ so that \lq tr' represents the
unique semifinite trace on ${\cal B}(\Omega_{(2)}(M))$ which extends the usual
trace on the trace class operators. Where no confusion can arise we simply write
$\tau$ for any of these traces.

Now we introduce the spectral density functions for these complexes. Let $d_j $
denote the restriction of $d$ to $j$-forms, with a similar convention for $
d^K_j$ and $d_j^{\widetilde K}$ acting on $C^j(K)$ and $C_{(2)}^j(\widetilde K)$
respectively. Let $\delta_j$ denote the Hilbert space adjoint of $d_j$ with
$\delta_j^{\widetilde K,\widetilde W}$ denoting the adjoint of the coboundary
$d^{\widetilde K}_j$ with respect to the inner product on $ C_{(2)}(\widetilde
K)$. Define the analytic Laplacian $\Delta_j$ to be equal to $d\delta + \delta
d$ acting on $\Omega_{(2)}(\widetilde M)$. We also define the combinatorial
Laplacian $\Delta_j^{\widetilde K,\widetilde W}$ to be equal to
$d^K\delta^{\widetilde K,\widetilde W} + \delta^{\widetilde K,\widetilde W} d^K$
acting on $C_{(2)}(\widetilde K)$. Regard $d_j$ as an operator on smooth
$j$-forms in $( {\rm ker}d_j)^\perp$. Then we may define the spectral resolution
for $ D_j\equiv\delta_jd_j=\int_0^\infty \lambda dE_\lambda^j$. By results of
Atiyah and Singer \cite{A,S} the spectral projections $E^j_\lambda$ are in the
commutant of the ${\cal U}$ action and have finite von Neumann trace. So we can
also define the spectral density function $F_j(\lambda)=\tau(E^j_\lambda) $.
Similar comments apply to the combinatorial operator $\delta_j^{ \widetilde
K,\widetilde W}d_j^{\widetilde K}$ so that we may similarly define its spectral
projections $E_\lambda ^{\widetilde K, \widetilde W, j}$ and the spectral
density function $F_j^{\widetilde K, \widetilde W}(\lambda)$.

\bigskip

{\em As in Dodziuk-Patodi \cite{DP}, we make the following assumptions on our
triangulations. Firstly, we consider only triangulations $K$ which are
rectilinear subdivisions of any fixed triangulation $K_0$. The next assumption
is that the fullness of $K$ (see  \cite{W} p. 125 for the definition) is bounded
away from zero.}

\bigskip

We will next state a version of an elementary lemma which is due to Novikov and
Shubin and whose proof can be found in the appendix of \cite{GS}.

Let $F:I\!\!R_+\rightarrow I\!\!R_{+}$ be a non-decreasing function satisfying
the following sub-exponential estimate, $ F(\lambda)=O(e^{\varepsilon \lambda})$
for all $\varepsilon>0$. Then the Laplace transform of $F$ exists, and is given
by \[ \theta(t)=\int^{\infty}_{0}e^{-\lambda t}dF(\lambda). \] Let
$\displaystyle\bar{b}=\lim_{\varepsilon\rightarrow 0^{+}}F(\varepsilon)=
F(0^{+})$. Then $\displaystyle\bar{b}=\lim_{t\rightarrow\infty} \theta(t)$.

\begin{lemma} : \cite{GS} $\;$Let $F$ and $\theta$ be as above, and let
$\alpha>0$. Then the following conditions are equivalent. \\[+7pt] $\!\!\!$(1).
There exists $C>0$ such that $$ C^{-1}\lambda^{\alpha}\leq
F(\lambda)-\bar{b}\leq C\lambda^{\alpha} $$ for all $\lambda$ small,
$\lambda>0$. \\[+7pt] $\!\!\!$(2).  There exists $C^{\prime}>0$ such that $$
C^{\prime}{^{-1}}t^{-\alpha}\leq\theta(t)-\bar{b}\leq C^{\prime}t^{-\alpha} $$
for all $t\gg 0$.

More generally, one has the following equality $$ \lim\inf_{\lambda\rightarrow
0^+}\Big\{ \frac{\log(F(\lambda)-\bar{b})} {\log \lambda} \Big\} =
\lim\inf_{t\rightarrow\infty}\Big\{ -\frac{\log(\theta(t) - \bar{b})}{\log t}
\Big\}. $$ \end{lemma}

In our context we need to consider a slightly different situation. Let $
G:I\!\!R\rightarrow I\!\!R_{+}$ be a non-decreasing function satisfying the
sub-exponential estimate $G(\lambda)=O(e^{\varepsilon\lambda})$ for all $
\varepsilon>0$, and $G(\lambda)=0$ for $\lambda<\lambda_{0}$ where $
\lambda_{0}>0$. Let $\psi(t)$ denote its Laplace transform, that is, \[
\psi(t)=\int^{\infty}_{\lambda_0}e^{-\lambda t}dG(\lambda)\;\;\mbox{ for }\;\;
t>0. \] Let
$\displaystyle\bar{a}=\lim_{\lambda\rightarrow\lambda_0^+}G(\lambda)=
\lim_{t\rightarrow\infty}\psi(t)$. Then the following lemma is an easy
consequence of Lemma 1.1.

\begin{lemma} : $\;$Let $G$ and $\psi$ be as above, and let $\alpha>0$. Then the
following conditions are equivalent. \\ $\!\!\!$(1). There exists $C>0$ such
that \[ C^{-1}(\lambda-\lambda_{0})^{\alpha}\leq G(\lambda)-\bar{a}\leq
C(\lambda- \lambda_{0})^{\alpha} \] for all $\lambda>\lambda_{0}$ and
$\lambda-\lambda_{0}$ small. \\ $\!\!\!$(2). There exists $C^{\prime}>0$ such
that \[ C^{\prime}{^{-1}}e^{-\lambda_{0}t} t^{-\alpha}\leq\psi(t)-\bar{a}\leq
C^{\prime}e^{-\lambda_{0}t} t^{-\alpha} \] for all $t\gg 0$.

More generally, one has the following equality $$
\lim\inf_{\lambda-\lambda_o\rightarrow 0^+}\Big\{
\frac{\log(G(\lambda)-\bar{a})} {\log(\lambda -\lambda_0)} \Big\} =
\lim\inf_{t\rightarrow\infty}\Big\{ -\frac{\log(\psi(t) - \bar{a})}{\log t } -
\frac{\lambda_0 t}{\log t} \Big\}. $$ \end{lemma}

In fact,    we shall  use below  the one-sided versions of the above Lemma (for
example, the RHS inequality in $(1)$ is equivalent to the RHS inequality in
$(2)$), and also the fact that the constants in $(1)$ and $(2)$ are connected:
if the situation depends on some parameter, $C$ and $C'$ go to $0$ or $+\infty$
at the same time.

\bigskip

\section{The Convergence of Spectral Functions}

In this section, we prove that the von Neumann spectral density  of the
combinatorial Laplacian converges to the von Neumann spectral density of the
analytic Laplacian as we successively refine our triangulations so that the mesh
goes to zero. We also prove the analogous theorem for the bottom of the spectra
of the combinatorial and analytic Laplacians. We use the notation $
H^s(\widetilde M)$ for the Sobolev space consisting of forms $\omega$ with $$
\int_{\widetilde M} \Vert (1+\Delta)^{s/2}\omega \Vert <\infty.\eqno (2.1) $$

For the readers' convenience we now recall (cf \cite{Do}) the definitions of the
de Rham and Whitney maps. For every continuous form $\omega$ of degree $p $ on
$\widetilde{M}$ we define a cochain $\int\omega\in C^{p}(\widetilde{K})$ by the
formula \[ \int\omega=\sum_\sigma\left(\int_{\sigma}\omega\right)\cdot\sigma. \]

\begin{lemma} : \cite{Do} Suppose $\omega\in H^{k}(\widetilde M)$,
$k>{\frac{n}{2}}+1$, has degree $p$. Define the de Rham map $\widetilde A\omega
= \int\omega$, then $ \widetilde A\omega$ is in $C_{(2)}(\widetilde K)$.
Moreover $\widetilde A:H^k(\widetilde M)\rightarrow C_{(2)}(\widetilde{K})$ is
bounded and $ \widetilde Ad=d^{\widetilde K}\widetilde A$. \end{lemma}

Let $\{U_{p}\}_{p\in K^{0}}$ be an open covering of $M$ by open stars of
vertices of $K$, that is, $U_{v}=$ st $(v)$ for $v\in K^0$. There is a partition
of unity on $M$ subordinate to this covering and which is given by the
barycentric coordinate functions. Both the covering and the partition of unity
can be lifted to $\widetilde{M}$ so that we obtain a partition of unity
$\{\mu_{v}\}_{v\in \widetilde{K}^{0}}$ indexed by vertices of $ \widetilde{K}$
with the following properties: $$ \begin{array}{l} \mbox{supp
}\mu_{v}\subset\mbox{supp\ st }v \\[+5pt]
\mu_{v}\cdot\gamma=\mu_{\gamma^{-1}v}\;. \end{array} \eqno(2.2) $$

Now let $\sigma=[v_{0},\ldots,v_{p}]$ be a simplex of $K$. To simplify the
notation we write $\mu_{i}=\mu_{v_i}$. Define the Whitney map

$$ \widetilde{W}\sigma=\left\{ \begin{array}{l} \displaystyle\mu_{0}\;\mbox{ if
}\;p=0 \\[+16pt] \displaystyle
p!\sum_{i=0}^{p}(-1)^{i}\mu_{i}d\mu_{0}\wedge\ldots\wedge d\mu_{i-1}\wedge
d\mu_{i+1}\wedge\ldots\wedge d\mu_{p}\;\mbox{ if } \; p>0. \end{array}
\right.\eqno(2.3) $$

Note that $\gamma^{*}\widetilde{W}\sigma=\widetilde{W}\gamma^{-1}\sigma$ by
(2.2). Another property of $\widetilde{W}\sigma$ which we need is $$ \mbox{supp
}\widetilde{W}\sigma=\mbox{ supp\ st }\sigma.\eqno(2.4) $$

Now for an arbitrary cochain $f=\Sigma f_{\sigma}c_\sigma$ we define $
\widetilde{W}f=\Sigma f_{\sigma}\widetilde{W}{\sigma}$. In view of (2.4), the
sum in (2.3) is locally finite and defines an $L^2$ form. For every $ f\in
C_{(2)}(\widetilde{K})$ the form $\widetilde Wf$ is in $\Omega_{(2)}(
\widetilde{M})$. Moreover $\;\widetilde W:C_{(2)}(\widetilde {K}) \rightarrow
\Omega_{(2)}(\widetilde M)$ is a bounded operator with the following properties.
$$ \begin{array}{ll} \mbox{(i) }\;\;\;\; & d\cdot
\widetilde{W}=\widetilde{W}\cdot d^{\widetilde K} \\[+5pt] \mbox{(ii) } &
\widetilde A\widetilde{W}=I\;\;\mbox{ on }\;\;C_{(2)}(\widetilde{K}) \\[+5pt]
\mbox{(iii) } & \gamma^{*}\widetilde{W}f=\widetilde{W}(f\cdot\gamma)\; \mbox{
for every }\;f\in C_{(2)} (\widetilde{K})\;\mbox{ and every }\;\gamma\in\Gamma.
\end{array} \eqno(2.5) $$

The first two properties are provided in \cite{W}. The third one is a
consequence of our choice of partition of unity.

Now we review some of the results of \cite{ES} and \cite{GS} which we need. Let
$H$ denote a Hilbert space with a trivial action of $\Gamma$ and let $
V\subseteq \ell^2(\Gamma)\otimes H$ be a closed $\Gamma$ invariant subspace. Let
$P_V$ denote the orthogonal projection onto $V$. Then we recall that $
\dim_\Gamma(V) = \tau(P_V)$. It can be seen that $\dim_\Gamma(V)$ is independent
of the choice of the $\Gamma$-invariant embedding of $V$ as a subspace of
$\ell^2(\Gamma)\otimes HS$, where $HS$ denotes a Hilbert space with a trivial
action of $\Gamma$. Let ${\cal S}^{j}_{\lambda}$ denote the set of $\Gamma$
invariant subspaces $L$ of $(\ker d_j)^{\bot}\subseteq\Omega_{(2)}^j(\widetilde
M)$ such that $$ ||d_j\omega||\leq\sqrt{\lambda}||\omega||\;\;\mbox{ for all
}\;\;\omega\in L \eqno (2.6) $$

and ${\cal S}^{j,\widetilde K,\widetilde W}_{\lambda}$ denote the set of $
\Gamma$ invariant subspaces $L$ of $(\ker d_j^{\widetilde K})^\bot\subseteq
C_{(2)}(\widetilde K)$ such that

$$ ||d_j^{\widetilde K}a||_{\widetilde W}\leq\sqrt{\lambda}||a||_{\widetilde
W}\;\;\;\forall\;a\in L. \eqno (2.7) $$

Then the respective spectral functions are also given by

$$ F_{j}(\lambda)\equiv \mbox{sup}\{\dim_{\Gamma}(L)\ \vert\ \ L\in {\cal S}
^{j}_{\lambda}\} \eqno (2.8) $$

and

$$ F_{j}^{\widetilde{K},\widetilde W}(\lambda)= \mbox{sup}\{\dim_{\Gamma}(L) \
\vert \ \ L\in{\cal S}^{j,\widetilde K,\widetilde W}_{\lambda}\} \eqno (2.9) $$

by the variational principle of \cite{ES} and \cite{GS}. This principle is the
analogue of the min-max principle in finite dimensions, and gives a variational
characterisation of the von Neumann spectrum of our Laplace operators. We
introduce the notation $\eta$ for the mesh of the triangulation $\widetilde K$
(which is of course the mesh for $K$).

Our first main result in this section is as follows,

\begin{proposition} : $\;$There is a constant $C_{1}$ independent of $\eta$ such
that if $ \displaystyle\mu<\left(\frac{1}{C_{1}\eta|\log\eta|}-1\right)^{2}$ and
$ \lambda\leq\mu$, then $$ F_{j}^{\widetilde{K},\widetilde{W}}(\lambda)\leq
F_{j}(D_\eta^\mu\lambda), $$ where  $\displaystyle
D_\eta^\mu=\Big(1-C_{1}\eta|\log\eta|(1+ \sqrt{\mu}\,) \Big)^{-2}\rightarrow 1$
as $\mu$ is fixed and $\eta\rightarrow 0$. \end{proposition}

The proof is contained in the following sequence of three lemmas. Let $P,
P^{\widetilde K}$ denote the projections onto $(\ker d_j)^\bot$ and $(\ker
d_j^{\widetilde K})^\bot$ respectively.

\begin{lemma} : $\;P\widetilde{W}$ is injective on $L$, where $L$ is a
$\Gamma$-invariant subspace of $(\ker d^{\widetilde{K}}_{j})^{\perp}$.
\end{lemma}

\noindent{\bf Proof}: $\;$Let $a\in L$ and $P\widetilde{W}a=0$, that is, $
\widetilde{W} a\in\ker d_j$. Then $$ \widetilde{W}d^{\widetilde{K}}_j a=d_j
\widetilde{W}a=0. $$

Hence $d^{\widetilde{K}}_j a=0$, since $\widetilde{W}$ is injective. But $
a\in(\ker d^{\widetilde{K}}_{j})^{\perp}$, so $a=0$.

An immediate consequence is that $$
\dim_{\Gamma}(L)=\dim_{\Gamma}(P\widetilde{W}(L)). $$ This is because
$P\widetilde W$ is injective on such subspaces and injective maps of ${\cal
U}$-modules into free modules preserve the $\Gamma$ dimension of the space.

\noindent Recall that $E^{\widetilde{K},\widetilde{W}}_{\mu}C^{j}_{(2)} (
\widetilde{K}) = \Big\{a\in C^{j}_{(2)}(\widetilde{K}) : ||d_j^{\widetilde{K} }
a||_{\widetilde{W}}^2\leq\mu||a||_{\widetilde{W}}^2\Big\}$.

The following is a highly non-trivial estimate in \cite{DP}.

\begin{lemma} : $\;$$\;$Let $a\in
E^{\widetilde{K},\widetilde{W}}_{\mu}C^{j}_{(2)} ( \widetilde{K})$, then one has
for $\eta$ sufficiently small $$
\|\widetilde{W}P^{\widetilde{K}}a-P\widetilde{W}a\|\leq C_{1}(1+\sqrt{\mu} \,)
\eta|\log\eta|\;\|a\|_{\widetilde{W}}. $$ \end{lemma}

\begin{lemma} : $\;$If $a\in E^{\widetilde{K},\widetilde{W}}_{\mu}C^{j}_{(2)}
(\widetilde{K })$ satisfies
$\|d^{\widetilde{K}}a\|_{\widetilde{W}}\leq\sqrt{\lambda}\;
\|a\|_{\widetilde{W}}$ and $a\in(\ker d^{\widetilde{K}}_j)^{\perp}$, then one
has $$
\|dP\widetilde{W}a\|\leq\Big(1-C_{1}\eta|\log\eta|(1+\sqrt{\mu}\,)\Big)^{-1}
\sqrt{\lambda}\;\|P\widetilde{W}a\|, $$

for $\eta$ sufficiently small with respect to $\mu$. \end{lemma}

\noindent{\bf Proof}: $\;$Since $d(1-P)=0$, we see that $$ \begin{array}{ll}
\|dP\widetilde{W}a\|= & \|\widetilde{W}d^{\widetilde{K}}a\|= \|d^{\widetilde{
K}} a \|_{\widetilde{W}} \\[+10pt] & \leq \sqrt{\lambda}\;\|a\|_{\widetilde{W}}
= \sqrt{\lambda}\;\|P^{ \widetilde{K}}a\|_{\widetilde{W}} \end{array} $$ since
by hypothesis $a = P^{\widetilde{K}}a$.

Lemma 2.4 and the equality $a = P^{\widetilde{K}}a$ imply that $$
\|P^{\widetilde{K}}a\|_{\widetilde{W}}\leq\;\Big(1-C_{1}\eta|\log\eta|(1+
\sqrt{\mu}\,)\Big)^{-1}\;\|P \widetilde{W}a\|. $$ This proves the lemma.

\noindent{\bf Proof of Proposition 2.2}: For fixed $\mu$ and for $\eta$
sufficiently small, let $\displaystyle D_\eta^\mu=\Big(1-C_{1}\eta|\log\eta|(1+
\sqrt{\mu}\,)\Big)^{-2}$.

It follows from Lemma 2.5 that if $L\in{\cal S}^{j,\widetilde{K},\widetilde{W
}}_\lambda$ and $\lambda\leq\mu$, then $P\widetilde{W}(L)\in
S^{j}_{D_\eta^\mu\lambda}\;$, and from Lemma 2.3 that $\dim_{\Gamma}(L)=\dim_{
\Gamma}(P\widetilde{W}(L))$.

The proposition is then merely a consequence of the definition of $F_j$ and $
F^{\widetilde{K},\widetilde{W}}_j$.

Our next main result of this section is as follows,

\begin{proposition} : There is a constant $C_2$ independent of $\eta$ such that
if $\mu<\big( C_2\eta\big)^{-2/s}-1$ and $\lambda\leq\mu$, then
$F_{j}(\lambda)\leq F_{j}^{ \widetilde{K},\widetilde W}(C_\eta^\mu\lambda)$,
where  $ \displaystyle
C_\eta^\mu=\left\{\frac{1+C_{2}\eta(1+\mu)^{\frac{s}{2}}}{1-C_{2} \eta(1+\mu)^{
\frac{s}{2}}}\right\}^2 \rightarrow 1$ as $\eta\rightarrow 0$. \end{proposition}

The proof relies on the following sequence of lemmas.

\begin{lemma} \cite{Do,DP} The inequality $||\omega-\widetilde W\widetilde
A\omega||\leq C_2\eta||\omega||_{s}$ holds for all $\omega \in
H^{s}(\widetilde{M})$ and $s>\frac{n}{2}+1\;$ where $C_2$ is independent of
$\eta$. \end{lemma}

Recall that $E_{\mu}\Omega^j_{(2)}(\widetilde M) = \Big\{\omega\in
\Omega^j_{(2)}(\widetilde M): ||d_j\omega||^2\leq\mu||\omega||^2\Big\}$.

The following lemma is taken from \cite{E}, and its proof is easy.

\begin{lemma} : For $\omega\in E_{\mu}\Omega^j_{(2)}(\widetilde M)$ and $s\geq
0$, $$ ||\omega||^{2}\leq||\omega||^{2}_{s}\leq(1+\mu)^{s}||\omega||^{2}. $$
\end{lemma}

It follows from this lemma that the de Rham map is well defined and bounded on
the range of $E_\mu$.

>From now on, we shall fix $s > {\frac{n}{2}} + 1$.

\begin{lemma} \cite{E}: For $\omega\in E_{\mu}\Omega^j_{(2)}(\widetilde M)\cap
(\ker d_j)^{\perp}$ and $\eta$ sufficiently small, one has $$
||\omega||\leq\frac{1}{1-C_2\eta(1+\mu)^{{\frac{s}{2}}}}\;||P^{\widetilde
K}\widetilde A\omega||_{\widetilde W}. $$ \end{lemma}

\noindent{\bf Proof:} The operator $$ B=(1-\widetilde W\widetilde
A)E_{\mu}:\Omega^{j}_{(2)}(\widetilde{M})
\rightarrow\Omega^{j}_{(2)}(\widetilde{M}) $$ is bounded for all $j$ with
$||B||<C_{2}\eta(1+\mu)^{\frac{s}{2}}$ (using Lemmas 2.7 and 2.8). Now because
$\omega=\widetilde W\widetilde A\omega+B\omega$, and also as $\widetilde W:\ker
d_j^{\widetilde K}\rightarrow\ker d_j$ and $\widetilde A:\ker d_j\rightarrow
\ker d_j^{\widetilde K}$, we conclude that $$ \begin{array}{lcl} ||\omega|| & =
& ||P\omega|| \leq ||P\widetilde W\widetilde A\omega||+||PB\omega|| \\[+5pt] &
\leq & ||P\widetilde WP^{\widetilde K}\widetilde A\omega||+ ||B\omega|| \\
[+5pt] & \leq & ||P^{\widetilde K}\widetilde A\omega||_{\widetilde
W}+C_{2}\eta(1+\mu)^{\frac{s}{2}}||\omega||. \end{array} $$ An immediate
consequence of this result is that the map $P^{\widetilde K}\widetilde A$ is
injective on $L$, where $L$ is a $\Gamma$ invariant subspace of
$E_{\mu}\Omega^j_{(2)}(\widetilde M)\cap (\ker d_j)^{\perp}$ as soon as $$ \eta<
{\frac{1}{C_{2}(1+\mu)^{\frac{s}{2}}}}. $$ Hence, for $\eta$ small enough with
respect to $\mu$, $$ \dim_{\Gamma}(L)=\dim_{\Gamma}(P^{\widetilde K}\widetilde
AL) $$ for all $\Gamma$-invariant subspaces $L$ in the range of $E_{\lambda}$
where $\lambda\leq\mu$. Again, this is because  injective maps of ${\cal
U}$-modules into free modules preserve the $\Gamma$ dimension of the space.

\begin{lemma} : Let $\omega\in E_{\mu}\Omega_{(2)}^j(\widetilde M)$ satisfy $
||d\omega||\leq\sqrt{\lambda}|| \omega||$ and $\omega\in (\ker d_j)^{\bot}$.
Then, for $\eta$ small enough, $$ ||d^{\widetilde K} P^{\widetilde K} \widetilde
A\omega||_{\widetilde
W}\leq\left\{\frac{1+C_2\eta(1+\mu)^{\frac{s}{2}}}{1-C_{2}\eta(1+\mu)^{\frac{
s}{2}}} \right\}\sqrt{\lambda}||P^{\widetilde K}\widetilde A\omega||_{\widetilde
W}. $$ \end{lemma}

\noindent{\bf Proof:} Notice first that since $d^{\widetilde K}(I -
P^{\widetilde K}) = 0$, one has $$ ||d^{\widetilde K} P^{\widetilde K}
\widetilde A\omega||_{\widetilde W} = ||d^{\widetilde K}\widetilde
A\omega||_{\widetilde W} = ||\widetilde Wd^{\widetilde K}\widetilde
A\omega||=||\widetilde W\widetilde Ad\omega||. $$ Now $d\omega\in
E_{\mu}\Omega_{(2)}^{j+1}(\widetilde M)$, since $d^2\omega = 0$. Therefore
according to Lemmas 2.8 and 2.9, $$ ||d\omega - \widetilde W\widetilde A
d\omega|| \leq C_2 \eta (1+\mu)^{\frac{s }{2}}||d \omega||. $$ Therefore using
the hypothesis $||d\omega||\leq\sqrt{\lambda}|| \omega||$, one has $$
||\widetilde W\widetilde A d\omega||\leq (1 + C_2 \eta (1+\mu)^{\frac{s}{2} })
\sqrt{\lambda}||\omega||. $$ The conclusion of the lemma follows from Lemma 2.9.

\noindent{\bf Proof of Proposition 2.6}: For fixed $\mu$ and for $\eta$
sufficiently small, let $\displaystyle
C_\eta^\mu=\left\{\frac{1+C_2\eta(1+\mu)^{\frac{s}{2}}}{1-C_2 \eta(1+\mu)^{
\frac{s}{2}}}\right\}^2$.

It follows from Lemma 2.10 that if $L\in S^{j}_{\lambda}$ and $\lambda\leq\mu $,
then $P^{\widetilde{K}}{\widetilde A} (L)\in {\cal S}^{j,\widetilde{K},
\widetilde{W}} _{C_\eta^\mu \lambda}\;$, and from Lemma 2.9 that $
\dim_{\Gamma}(L)=\dim_{\Gamma}(P^{\widetilde K}\widetilde AL)$.

The proposition is then merely a consequence of the definition of $F_j$ and $
F^{\widetilde{K},\widetilde{W}}_j$.

We will now prove that ${\kappa}^{\widetilde{K}, \widetilde W}_{0,j} $, the
bottom of the spectrum of $\delta_j^{\widetilde K, \widetilde W}d_j^{\widetilde
K}$ acting on $\delta_j^{\widetilde K,\widetilde{W}
}C^{j+1}_{(2)}(\widetilde{K})$,  converges to $\kappa_{0,j}$, the bottom of the
spectrum of $\delta_j d_j$ acting on $\delta_j
\Omega^{j+1}_{(2)}(\widetilde{M})$,  as the mesh of the triangulation $K$ goes
to zero.

First notice that $\overline{\delta_j^{\widetilde{K},\widetilde{W}
}C^{j+1}_{(2)}(\widetilde{K})} =(\ker d^{\widetilde{K}}_j)^{\perp}=(\ker
\delta_j^{\widetilde{W},\widetilde{K}} d_j^{\widetilde{K}})^{\perp}$, and that
$\overline{\delta_j \Omega^{j+1}_{(2)}(\widetilde{M})} =(\ker
d_j)^{\perp}=(\ker\delta_j d_j)^{\perp}$.

Therefore $$ \begin{array}{llcl} & \kappa^{\widetilde{K},\widetilde W}_{0,j} & =
& \displaystyle\inf \left\{ \frac{\|d^{\widetilde{K}}{a} \|^{2}_{\widetilde
W}}{\|a\|^{2}_{\widetilde W}} \;:\;\;a\in (\ker
d^{\widetilde{K}}_j)^{\perp}\right\} \\[+16pt] \mbox{and }\;\; & \kappa_{0,j} &
= & \displaystyle\inf\left\{\frac{ \|d\omega\|^{2}}
{\|\omega\|^{2}}\;:\;\;\omega\in (\ker d_j)^{\perp}\right\}. \end{array} $$

\begin{theorem} : The following inequality holds
$$
D_\eta^{-1}\kappa_{0,j}\leq\kappa^{\widetilde{K}, \widetilde W}_{0,j}
\leq C_\eta\kappa_{0,j}, $$ where $D_\eta$ and $C_\eta$ tend to $1^+$ as
$\eta$ goes
to zero.

In particular, $\;\kappa^{\widetilde{K}, \widetilde W}_{0,j}$ converges to $
\kappa_{0,j}$ as the mesh $\eta$ of the triangulation $K$ goes to zero.
\end{theorem}

\noindent{\bf Proof}: $\;$ Take $\varepsilon>0$, and let $\omega\in (\ker
d_j)^{\perp}$ be such that $$||d\omega||\leq(\kappa_{0,j}+\varepsilon)^{1/2}||
\omega||.$$ According to Lemma 2.10, for $\eta$ small enough, $$ ||d^{\widetilde
K} P^{\widetilde K} \widetilde A\omega||_{\widetilde
W}\leq\left\{\frac{1+C_2\eta(1+\kappa_{0,j}+\varepsilon)^{\frac{s}{2}}}{1-C_
{2}\eta(1+\kappa_{0,j}+\varepsilon)^{\frac{ s}{2}}}
\right\}\sqrt{\kappa_{0,j}+\varepsilon}||P^{\widetilde K}\widetilde
A\omega||_{\widetilde W}. $$ Therefore $$\kappa^{\widetilde{K}, \widetilde
W}_{0,j}\leq
\left\{\frac{1+C_2\eta(1+\kappa_{0,j}+\varepsilon)^{\frac{s}{2}}}
{1-C_{2}\eta(1+\kappa_{0,j}+\varepsilon)^{\frac{s}{2}}}
\right\}^2(\kappa_{0,j}+\varepsilon). $$ Since $\varepsilon$ is arbitrary, this
gives $$\kappa^{\widetilde{K}, \widetilde W}_{0,j}\leq
\left\{\frac{1+C_2\eta(1+\kappa_{0,j})^{\frac{s}{2}}}{1-C_{2}\eta(1+\kappa_{
0,j})^{\frac{ s}{2}}} \right\}^2\kappa_{0,j},$$ and one  sets $$
C_\eta=\left\{\frac{1+C_2\eta(1+\kappa_{0,j})
^{\frac{s}{2}}}{1-C_{2}\eta(1+\kappa_{0,j})^{\frac{s}{2}}}\right\}^2. $$ Take
$\varepsilon>0$, and let $a\in (\ker d^{\widetilde{K}}_j)^{\perp}$ be such that
$$\|d^{\widetilde{K}}a\|_{\widetilde{W}}\leq (\kappa^{\widetilde{K}, \widetilde
W}_{0,j}+ \varepsilon)^{1/2}\; \|a\|_{\widetilde{W}}.$$ For $\eta$ small enough,
$\Big(1-C_{1}\eta|\log\eta|(1+(\kappa^{\widetilde{K}, \widetilde W}_{0,j}+
\varepsilon)^{1/2}\,)\Big)$ is positive, since by the inequality just proven
$\kappa^{\widetilde{K}, \widetilde W}_{0,j}$ is bounded when $\eta$ goes to
zero. According to Lemma 2.5,
$$\|dP\widetilde{W}a\|\leq\Big(1-C_{1}\eta|\log\eta|(1+(\kappa^{\widetilde{K },
\widetilde W}_{0,j}+ \varepsilon)^{1/2}\,)\Big)^{-1} (\kappa^{\widetilde{K},
\widetilde W}_{0,j}+ \varepsilon)^{1/2}\;\|P\widetilde{W}a\|, $$ hence $$
\kappa_{0,j} \leq\Big(1-C_{1}\eta|\log\eta|(1+\sqrt{\kappa^{\widetilde{K},
\widetilde W}_{0,j}}\,)\Big)^{-2} \kappa^{\widetilde{K}, \widetilde W}_{0,j}. $$
Now, if $\kappa^{\widetilde{K}, \widetilde W}_{0,j}$ is bounded from above by
$K$ as $\eta$ goes to zero, one sets
$$D_\eta=\Big(1-C_{1}\eta|\log\eta|(1+\sqrt{K}\,)\Big)^{-2}.$$

According to \cite{DP}, p.11 (see also \cite{GS} p.385), the spectrum of
$\delta_j d_j$ acting on $\delta_j \Omega^{j+1}_{(2)}(\widetilde{M})$ coincides
with the spectrum of $d_j\delta_j$ acting on $d_j
\Omega^{j}_{(2)}(\widetilde{M})$. Therefore one has $$ \lambda_{0,j} = minimum\{
\kappa_{0,j-1}, \kappa_{0,j}\}, $$ where $\lambda_{0,j}$ denotes the bottom of
the spectrum of the Laplacian $\Delta_{j}$ acting on $\Omega^{j}_{(2)}(
\widetilde{M})$. Similar remarks apply in the combinatorial situation, and we
can therefore state

\begin{corollary} : $\;$ The bottom of the spectrum of the combinatorial
Laplacian $\Delta^{ \widetilde{K},\widetilde{W}}_{j}$ acting on $C^{j}_{(2)}
(\widetilde{K})$, $ \lambda^{\widetilde{K},\widetilde W}_{0,j}$, converges to
the bottom of the spectrum of the Laplacian $\Delta_{j}$ acting on
$\Omega^{j}_{(2)}( \widetilde{M})$, $\lambda_{0,j}$. More precisely, the
following inequality holds $$
{D_\eta}^{-1}\lambda_{0,j}\leq\lambda^{\widetilde{K}, \widetilde W}_{0,j}\leq
C_\eta\lambda_{0,j}. $$ where $D_\eta$ and $C_\eta$ tend to $1^+$ as $\eta$ goes
to zero. \end{corollary}

\bigskip

\section{Estimates for L$^2$ Theta Functions.}

In this section, we prove some basic estimates for the $L^2$ theta functions
which we will need later on in the paper.

Note firstly that if $\bar N_j = \tau(E_{\Delta_j})$, where $\Delta_j =\int
\lambda dE_{\Delta_j}(\lambda)$ is the spectral resolution of the Laplacian $
\Delta_j$, then $\bar N_j=F_{j-1}+F_j + b^j_{(2)}(\widetilde M)$ (cf. \cite
{GS}). One sets $F_{-1} = F_{n+1} = 0$. Let ${\bar N_j} = N_j +
b^j_{(2)}(\widetilde M)$, that is, $N_j=F_{j-1}+F_j$. Recall the definition of
the theta function for the Laplacian on $L^2\;j$-forms: $$
\theta_{j}(t)=\tau(e^{-t\Delta_{j}}) - b^{j}_{(2)}(\widetilde
M)=t\int^{\infty}_{0}e^{-t\lambda}N_{j}(\lambda)d\lambda. $$ Since $$
\frac{d}{dt}\;\theta_{j}(t)=-\tau(\Delta_{j}e^{-t\Delta_{j}})<0, $$ $\theta_{j}$
is a decreasing function.

We use a similar notation for the theta function of the combinatorial Laplacian,
that is, $$ \theta_{j}^{\widetilde{K},\widetilde W}(t)= t\int_0^\infty
e^{-t\lambda}N_{j}^{\widetilde{K},\widetilde W} (\lambda)d\lambda. $$ For the
needs of section 4, we consider the truncated theta functions $$
\theta_{j,\nu}(t)=t\int^{\infty}_{\lambda_{0,j}+\nu_{\eta}}e^{-\lambda t}N_{j}
(\lambda)d\lambda. $$ and $$
\theta^{\widetilde{K},\widetilde{W}}_{j,\nu}(t)=t\int^{\infty}_{ \lambda_{0,j}+
\nu_{\eta}}e^{-\lambda t}N^{\widetilde{K},\widetilde{W} }_{j}(\lambda)d\lambda
$$ where $\nu_{\eta}$ is non-negative and bounded. Our aim is to investigate the
convergence of the former to the latter as the mesh goes to zero. Note that
$\theta_{j}^{\widetilde{K},\widetilde W}$, $\theta_{j,\nu}$,
$\theta^{\widetilde{K},\widetilde{W}}_{j,\nu}$ are decreasing as well.

\begin{proposition} There is a $t_0=t_0(\eta)$ such that if $t\geq t_{0}>0$,
then $$ \begin{array}{l} (1).\;\;\;\; e^{\lambda_{0,j} {D_\eta}^{-1}
t}\theta_{j,\nu}(C'_\eta t)\leq e^{\lambda^{\widetilde{K},\widetilde{W}}_{0,j}
t} \theta_{j,\nu}^{\widetilde{ K},\widetilde W}(t)+\varepsilon(\eta) \\ \\
(2).\;\;\;\; e^{\lambda_{0,j} {D_\eta}^{-1} t}\theta_{j}(C'_\eta t)\leq
e^{\lambda^{\widetilde{K},\widetilde{W}}_{0,j} t} \theta_{j}^{\widetilde{K}
,\widetilde W}(t)+\varepsilon(\eta) \\ \\ (3).\;\;\;\; \theta_{j}(C'_\eta t)\leq
\theta_{j}^{\widetilde{K},\widetilde W}( t)+\varepsilon(\eta) \end{array} $$
where $t_0\rightarrow 0, \ \ \varepsilon\rightarrow 0$ as $\eta\rightarrow 0$,
$D_\eta$ is as in Proposition 2.12 and $C'_\eta$ tends to $1^+$ as $\eta$ goes
to zero. \end{proposition}

\noindent{\bf Proof:} Choose $\mu=\eta^{-\alpha}$ , $0<\alpha<2/s$. Then, for
$\eta$ small enough, $C_\eta^\mu$ exists and tends to one as $\eta$ goes to
zero. Define $C'_\eta=C_\eta^\mu$. Write $$ e^{\lambda_{0,j} {D_\eta}^{-1}
t}\theta_{j,\nu}(C'_\eta t)=C'_\eta t e^{\lambda_{0,j} {D_\eta}^{-1} t}
\int^{\mu}_{\lambda_{0,j}+\nu_{\eta}}e^{-tC'_\eta\lambda}N_{j}(\lambda)d\lambda+
C'_\eta t e^{\lambda_{0,j} {D_\eta}^{-1} t}\int^{\infty}_{\mu}
e^{-tC'_\eta\lambda}N_{j}(\lambda)d\lambda. $$ By Propositions 2.6 and 2.12, $$
\begin{array}{lcl} \displaystyle C'_\eta t e^{\lambda_{0,j} {D_\eta}^{-1} t}
\int^{\mu}_{\lambda_{0,j}+ \nu_{\eta}}e^{-t
C'_\eta\lambda}N_{j}(\lambda)d\lambda & \leq & \displaystyle C'_\eta t
e^{\lambda^{\widetilde{K},\widetilde{W}}_{0,j} t}
\int^{\mu}_{\lambda_{0,j}+\nu_{\eta}}e^{-t C'_\eta\lambda}N_{j}^{\widetilde{K}
,\widetilde W}(C'_\eta \lambda)d \lambda \\[+12pt] & \leq & \displaystyle
e^{\lambda^{\widetilde{K},\widetilde{W}}_{0,j}
t}\theta^{\widetilde{K},\widetilde W}_{j,\nu}( t). \end{array} $$ We now proceed
to estimate the error \[ \varepsilon(\eta,t)\equiv C'_\eta
t\int^{\infty}_{\mu}e^{-t(C'_\eta\lambda - {D_\eta}^{-1}
\lambda_{0,j})}N_{j}(\lambda) d\lambda. \] Recall that the small time
asymptotics of the heat kernel on differential forms on the covering space is
given for $n=\ $dim$\ M$ by \[ \tau(e^{-t\Delta_j})=c_jt^{-n/2}+R(t) \] where
$\lim_{t\rightarrow 0}t^{n/2}R(t) =0$ as in \cite{R}. Hence by the Tauberian
theorem, the function $N_j(\lambda)$ has the large $\lambda$ asymptotic
expansion \[ N_j(\lambda)=c_j\lambda^{n/2}+f_j(\lambda) \] where
$\lim_{\lambda\rightarrow \infty}\lambda^{-n/2}f_j(\lambda ) =0$. It follows
that, given $\epsilon>0$, there is a $\Lambda_j>0$ such that for all $
\lambda>\Lambda_j$ \[ N_{j}(\lambda)\leq (c_{j}+\epsilon)\lambda^{n/2} . \]
Hence, if $t\geq t_0={\frac{1}{\sqrt{C'_\eta\mu - {D_\eta}^{-1}\lambda_{0,j}
}}}$, then $\varepsilon(\eta, t)\leq\varepsilon(\eta,t_0)$. Now, $\mu$ tends to
infinity as $\eta$ goes to zero. Thus, for $\eta$ small enough, \[
\begin{array}{lcl} \displaystyle\varepsilon(\eta,t_{0}) & \leq & \displaystyle
C'_\eta t_0(c_{j}+\epsilon)\int^{\infty}_ {\mu}e^{-t_0(C'_\eta\lambda -
{D_\eta}^{-1}\lambda_{0,j})} \lambda^{n/2}d\lambda \\[+12pt] & \leq &
\displaystyle C'_\eta t_0(c_{j}+\epsilon)e^{-\frac{t_0}{2}(C'_\eta\mu -
{D_\eta}^{-1}\lambda_{0,j})} \;\int^{\infty}_{\mu}e^{- \frac{t_0}{2}
(C'_\eta\lambda - {D_\eta}^{-1}\lambda_{0,j})}\lambda^{n/2}d\lambda \\ [+12pt] &
\leq & \displaystyle c_{j}^{\prime}(C'_\eta)^{-n/2}{(C'_\eta\mu -
{D_\eta}^{-1}\lambda_{0,j})} ^{\frac{n}{4}} e^{-\frac{\sqrt{C'_\eta\mu -
{D_\eta}^{-1}\lambda_{0,j}}}{2} }\rightarrow 0\;\;\mbox{ as }\;\;\eta\rightarrow
0 \end{array} \] where \[
c^{\prime}_{j}=(c_{j}+\epsilon)2^{\frac{n}{2}+1}\int_0^\infty e^{-y} (y+
\lambda_{0,j})^{\frac{n}{2}} dy. \] The proof above also works for parts (2) and
(3).

As a converse to the previous result we establish the following inequality.

\begin{proposition} There is a $t_0=t_0(\eta)$ such that if $t\geq t_{0}>0$,
then \[ \begin{array}{l} (1).\;\;\;\;
e^{\lambda^{\widetilde{K},\widetilde{W}}_{0,j}t}\theta_{j,\nu}^{
\widetilde{K},\widetilde W}(t)\leq e^{\lambda_{0,j}C_\eta
t}\theta_{j,\nu}(D'^{-1}_\eta t)+\varepsilon(\eta) \\ \\ (2).\;\;\;\;
e^{\lambda^{\widetilde{K},\widetilde{W}}_{0,j}t}\theta_{j}^{
\widetilde{K},\widetilde W}(t)\leq e^{\lambda_{0,j}C_\eta
t}\theta_{j}(D'^{-1}_\eta t)+\varepsilon(\eta) \\ \\
(3).\;\;\;\;\theta_{j}^{\widetilde{K},\widetilde W}(t)\leq
\theta_{j}(D'^{-1}_\eta t)+\varepsilon(\eta) \end{array} \] where
$t_0\rightarrow 0, \ \ \varepsilon\rightarrow 0$ as $\eta$ goes to zero,
$C_\eta$ is as in Proposition 2.12 and $D'_\eta$ tends to $1^+$ as $\eta$ goes
to zero. \end{proposition}

\noindent{\bf Proof:} Choose for example $\mu=\eta^{-\alpha}$, $0<\alpha<1$.
Then, for $\eta$ small enough, $D_\eta^\mu$ exists and tends to one as $\eta$
goes to zero. Define $D'_\eta=D_\eta^\mu$. Write \[
e^{\lambda^{\widetilde{K},\widetilde{W}}_{0,j}t}\theta_{j,\nu}^ {\widetilde{K
},\widetilde W}(t)=e^{\lambda^{\widetilde{K},\widetilde{W}}_{0,j}t}
t\int^{\mu}_{\lambda_{0,j}+\nu_{\eta}}e^{-t\lambda} N^{\widetilde{K} ,\widetilde
W}_{j}(\lambda)d\lambda+ e^{\lambda^{\widetilde{K},\widetilde{W}
}_{0,j}t}t\int^{\infty}_{\mu} e^{-t\lambda}N^{\widetilde{K},\widetilde
W}_{j}(\lambda)d\lambda. \] By Propositions 2.2 and 2.12, \[ \begin{array}{lcl}
\displaystyle e^{\lambda^{\widetilde{K},\widetilde{W}}_{0,j}t}
t\int^{\mu}_{\lambda_{0,j}+\nu_{\eta}}e^{-t\lambda} N^{\widetilde{K} ,\widetilde
W}_{j}(\lambda)d\lambda & \leq & \displaystyle e^{\lambda^{
\widetilde{K},\widetilde{W}}_{0,j}t}
t\int^{\mu}_{\lambda_{0,j}+\nu_{\eta}}e^{-t\lambda}N_{j}(D'_\eta \lambda)d
\lambda \\[+12pt] & \leq & \displaystyle e^{\lambda_{0,j}C_\eta
t}\theta_{j,\nu}(D'^{-1}_\eta t). \end{array} \] We now proceed to estimate the
error \[ \varepsilon(\eta,t)\equiv t\int^{\infty}_{\mu}e^{-t(\lambda- \lambda^{
\widetilde{K},\widetilde{W}}_{0,j})}N^{\widetilde{K},\widetilde W}_{j}(\lambda)
d\lambda. \] By the fullness assumption on the triangulation $\widetilde K$,
$\dim_\Gamma (C^{j}_{(2)}(\widetilde{K})) = \dim(C^{j}({K})) $ grows like
$\eta^{-n}$ as $ \eta\rightarrow 0$, see \cite{DP}. One then sees that the
spectral density function of the combinatorial Laplacian on $L^2\;j$-cochains
$N^{\widetilde{K },\widetilde W}_j(\lambda)$ must be bounded, \[
N^{\widetilde{K},\widetilde W}_j(\lambda) \leq c_j \eta^{-n}. \] Hence, if
$t\geq t_0={\frac{1}{{\sqrt{\mu -\lambda^{\widetilde{K}, \widetilde{W}}_{0,j}}
}}}$, then $\varepsilon(\eta, t)\leq\varepsilon(\eta,t_0)$ and \[
\begin{array}{lcl} \displaystyle\varepsilon(\eta,t_{0}) & \leq & \displaystyle
t_0 c_{j}\eta^{-n}\int^{\infty}_ {\mu}e^{-t_0(\lambda - \lambda^{\widetilde{K},
\widetilde{W}}_{0,j})} d\lambda \\[+12pt] & = & \displaystyle c_{j} \eta^{-n}
e^{-t_0 (\mu -\lambda^{\widetilde{K}, \widetilde{W}}_{0,j})} \\[+12pt] & = &
\displaystyle c_{j} \eta^{-n} e^{-{\sqrt{\mu - \lambda^{\widetilde{K
},\widetilde{W}}_{0,j}}}} \rightarrow 0\;\;\mbox{ as }\;\;\eta\rightarrow 0.
\end{array} \] The proof above also works for parts (2) and (3).

Finally, since  $\theta_j$, $\theta_{j}^{\widetilde{K},\widetilde W}$,
$\theta_{j,\nu}$, and $\theta^{\widetilde{K},\widetilde{W}}_{j,\nu}$ are
decreasing, one can replace $C_\eta$ and $C'_\eta$ by $\sup(C_\eta,C'_\eta)$,
$D_\eta$ and $D'_\eta$ by $\sup(D_\eta,D'_\eta)$ in 2.12, 3.1, and 3.2; this
means that from now on we can identify $C_\eta$ and $C'_\eta$, $D_\eta$ and
$D'_\eta$, as long as we only use the fact that they tend to one from above as
$\eta$ goes to zero.

\bigskip

\section{The Main Approximation Theorem}

In this section, we prove our main approximation theorem, Theorem 4.1. This is
proved using the approximation theorems for the truncated $L^2$ theta functions
of the previous section, as well as the approximation theorem for the bottom of
the spectra of Laplacians in section 2. We will assume that $\lambda_{0,j}>0$ in
this section. The case when  $\lambda_{0,j}=0$ is proved in section 7.

Recall the definition of the $\beta$ invariant, \[ \begin{array}{lcl}
\mbox{{$\beta$}}_{j}(M,g) = \displaystyle  \sup\{\beta\in
I\!\!R:e^{\lambda_{0,j}t}\theta_{j,\Delta}(t)\;\;\mbox{ is }
\;\;O(t^{-\beta})\;\mbox{ as }\;t\rightarrow\infty\}\in[0,\infty] &  & \\[+16pt]
=\displaystyle\liminf_{t\rightarrow \infty}\left\{-\frac{\log\theta_{j}(t)}{\log
t}-\frac{\lambda_{0,j}t}{\log t} \right\} \end{array} \] Then our main theorem
is

\begin{theorem} : $\;$ Assume that $\beta_j(M, g) >0$. Then in the limit as the
mesh $\eta$ of the triangulation $K$ goes to zero,  $e^{\lambda^{\widetilde{K},
\widetilde{W}}_ {0,j}t}\theta^{\widetilde{K},\widetilde{W}}_{j}(t)$ converges
uniformly to $e^{\lambda_{0,j}t}\theta_{j}(t)$ on $[t_0,+\infty[$ for any
$t_0>0$. \end{theorem}

The proof of this theorem is an immediate consequence of the following three
lemmas. For the rest of this section, we assume that $\beta_j(M, g) >0$. First
choose $\nu_{\eta}$ satisfying the following conditions:

(1) $\;\nu_{\eta}>0$ and $\nu_{\eta}\rightarrow 0$ as $\eta\rightarrow 0$

(2) $\;\lambda_{0,j}^{\widetilde{K},\widetilde{W}}<\lambda_{0,j}+\nu_{\eta}$

(3) $\;D^{-1}_{\eta}\lambda-\lambda_{0,j}C_{\eta}>0$ for $\lambda\geq
\lambda_{0,j}+\nu_{\eta}\;$.

Here $D_\eta,C_\eta\rightarrow 1^+$ as $\eta\rightarrow 0$, and satisfy 2.12,
3.1 and 3.2 (see the remark at the end of section 3). Then we have

\begin{lemma} : As $\eta$ goes to zero, $e^{\lambda_{0,j}t}\theta_{j,\nu}(t) -
e^{\lambda^{\widetilde{K}, \widetilde{W}}_
{0,j}t}\theta^{\widetilde{K},\widetilde{W}}_{j,\nu}(t)$ converges uniformly to
$0$ on $[t_1,+\infty[$ for any $t_1>0$. \end{lemma}

\noindent{\bf Proof}: By Proposition 3.2 and Theorem 2.12, one has for $t_0 =
t_0(\eta)$ and  $t\geq t_0$ \[
e^{\lambda^{\widetilde{K},\widetilde{W}}_{0,j}t}\theta^{\widetilde{K},
\widetilde{W}}_ {j,\nu}(t)-e^{\lambda_{0,j}t}\theta_{j,\nu}(t) \\[+10pt] \leq\;
e^{\lambda_{0,j}C_
{\eta}t}\theta_{j,\nu}(D^{-1}_{\eta}t)-e^{\lambda_{0,j}t}\theta_{j,\nu}(t)+
\varepsilon(\eta). \] Also by Proposition 3.1 and Theorem 2.11, one has for
$\eta$ small enough and $t_0 = t_0(\eta)$ and for $t\geq t_0$ \[
e^{\lambda_{0,j}t}\theta_{j,\nu}(t)-e^{\lambda^{\widetilde{K},\widetilde{W} }_
{0,j}t}\theta^{\widetilde{K},\widetilde{W}}_{j,\nu}(t) \leq
e^{\lambda_{0,j}t}\theta_{j,\nu}(t)-e^{\lambda_{0,j}D^{-1}_{\eta}t}
\theta_{j,\nu}(C_\eta t)+\varepsilon(\eta). \] So we deduce that for $t\geq t_0$
\[ \left|e^{\lambda_{0,j}t}\theta_{j,\nu}(t)-e^{\lambda^{\widetilde{K},
\widetilde{W}}_ {0,j}t}\theta^{\widetilde{K},\widetilde{W}}_{j,\nu}(t)\right|
\leq\;e^{\lambda_{0,j}C_{\eta}t}\theta_{j,\nu}(D^{-1}_{\eta}t)-e^{ \lambda_
{0,j}D^{-1}_{\eta}t}\theta_{j,\nu}(C_\eta t)+2\varepsilon(\eta). \] Let us
estimate \[ \begin{array}{l} I(t,\eta) =
e^{\lambda_{0,j}C_{\eta}t}\theta_{j,\nu}(D^{-1}_{\eta}t)-e^{ \lambda_{0,j}
D^{-1}_{\eta}t}\theta_{j,\nu}(C_\eta t) \\[+16pt]
\displaystyle=t\int^{\infty}_{\lambda_{0,j}+\nu}e^{-(\lambda D^{-1}_{\eta}-
\lambda_{0,j}C_{\eta})t}N_{j}(\lambda)d\lambda-t\int^{\infty}_{ \lambda_{0,j}
+\nu}e^{-(\lambda C_\eta-\lambda_{0,j}D^{-1}_{\eta})t}N_{j}(\lambda)d\lambda.
\end{array} \] Let \[ \begin{array}{llcl} & \alpha(\lambda,\eta) & = & \lambda
D^{-1}_{\eta}-\lambda_{0,j}C_{\eta} \\[+10pt] \mbox{and }\;\;\; &
\beta(\lambda,\eta) & = &\lambda C_\eta- \lambda_{0,j}D^{-1}_ {\eta}.
\end{array} \] First fix $\mu_0$ so that $t_1>1/(\mu_0-\lambda_{0,j})$. Then \[
\begin{array}{ll} & I(t,\eta)=I_1(t,\eta)+I_2(t,\eta) \\[+16pt] \mbox{where
}\;\;\; & \displaystyle I_1(t,\eta)=\int^{\mu_0}_{\lambda_{0,j}+\nu} t
\big(e^{-\alpha(\lambda, \eta)t}-e^{-\beta(\lambda,\eta)t}\big)
N_{j}(\lambda)d\lambda \\[+16pt] \mbox{and }\;\;\; & \displaystyle
I_2(t,\eta)=\int^{\infty}_{\mu_0}t\big( e^{-\alpha
(\lambda,\eta)t}-e^{-\beta(\lambda,\eta)t}\big) N_{j}(\lambda)d\lambda.
\end{array} \]

Let us first examine $I_2(t,\eta)$. For $\eta$ near zero,  the function
$t\rightarrow t(e^{-\alpha(\lambda,\eta) t}- e^{-\beta(\lambda,\eta) t})$ is
monotonic decreasing for $t >1/\alpha$. Now, for $\eta$ small enough,
$1/\alpha<t_1$, and for $t\geq t_1$, the integrand in $I_2(t,\eta)$ is smaller
than $t_1(e^{-\alpha(\lambda,\eta)t_1}-e^{- \beta(\lambda,\eta)t_1})$. Thus
$I_2(t,\eta)< I_2(t_1,\eta)$, $\; t\geq t_1$ and we have to show the right hand
side can be made small.

Now consider $r(\eta)= t_1( e^{-\alpha(\lambda,\eta)
t_1}-e^{-\beta(\lambda,\eta) t_1})$ as a function of $\eta$. As $
\eta\rightarrow 0$, $\alpha(\lambda,\eta)$ increases to $\lambda- \lambda_{0,j}$
while $\beta(\lambda,\eta)$ decreases to the same value. Thus  $r(\eta)$
decreases to zero as  $\eta$ goes to zero. By the dominated convergence theorem,
$I_2(t_1,\eta)$ goes to zero. We conclude that $I_2(t, \eta)\rightarrow 0$ {\em
uniformly} in $t$, $t\geq t_1$ as $\eta\rightarrow 0$.

Next we estimate $I_1(t,\eta)$. Consider \[
f(t)=te^{-\alpha(\lambda,\eta)t}\big(1-e^{-(\beta(\lambda,\eta)-\alpha( \lambda,
\eta))t}\big) \] for $t\in\big[0,\frac{2}{\alpha}\big]$. Then
$\max_{t\in\big[0,+\infty\big]}f(t)=\max_{t\in\big[0,\frac{2}{\alpha}\big]}f
(t)$ and \[ \begin{array}{lcl}
\displaystyle\max_{t\in\big[0,\frac{2}{\alpha}\big]}f(t) & \leq &
\displaystyle\max_{t\in\big[0,\frac{2}{\alpha}\big]}te^{-\alpha(\lambda,
\eta)t}\max_{t\in\big[0,\frac{2}{\alpha}\big]} (1-e^{-(\beta-\alpha)t}) \\
[+16pt] & \leq & \displaystyle\frac{1}{e\cdot\alpha}\;\frac{2}{\alpha}
\;(\beta-\alpha). \end{array} \] So \[ I_1(t,\eta) \leq
\displaystyle\frac{2}{e}\int^{\mu_0}_{\lambda_{0,j}+\nu} \;
\frac{(\beta(\lambda,\eta)-\alpha(\lambda,\eta))}{\alpha(\lambda,\eta)^{2}}
\;N_{j}(\lambda)d\lambda. \] Now
\[\beta(\lambda,\eta)-\alpha(\lambda,\eta)=(C_\eta - {D_\eta}^{-1})(\lambda +
\lambda_{0,j}).\] Moreover, since $\beta_j(M,g)>0$, Lemma 1.2 tells us that
\[N_j(\lambda)\leq C(\lambda-\lambda_{0,j})^{\beta},\] for $\beta\in ]0,1[$ and
$\lambda\leq \lambda_{0,j}+1$. We first estimate \[
\frac{2}{e}\int^{\lambda_{0,j}+1}_{\lambda_{0,j}+\nu} \frac{
(\beta(\lambda,\eta)-\alpha(\lambda,\eta))}{\alpha(\lambda,\eta)^{2}}
N_{j}(\lambda)d\lambda \leq \frac{\bar C}{1-\beta} (2\lambda_{0,j} +1)(C_\eta -
{D_\eta}^{-1})\Big( \nu_\eta^{\beta-1} -1\Big). \] Now assumption $(3)$ implies
that $\frac{(C_\eta -{D_\eta}^{-1})}{\nu_\eta}$ is bounded from above. Using
this we see that \[
\displaystyle\frac{2}{e}\int^{\lambda_{0,j}+1}_{\lambda_{0,j}+\nu} \frac{
(\beta(\lambda,\eta)-\alpha(\lambda,\eta))}{\alpha(\lambda,\eta)^{2}}
N_{j}(\lambda)d\lambda \\[+16pt] \displaystyle\leq\frac{\bar C}{1-\beta}
(2\lambda_{0,j} +1)\Big( C \nu_\eta^{\beta} -(C_\eta -{D_\eta}^{-1})\Big). \]
This is a uniform estimate in $t$, $t\geq t_1$, and the right hand side goes to
zero as $\eta\rightarrow 0$. We next estimate \[
\frac{2}{e}\int^{\mu_0}_{\lambda_{0,j}+1} \frac{
(\beta(\lambda,\eta)-\alpha(\lambda,\eta))}{\alpha(\lambda,\eta)^{2}}
N_{j}(\lambda)d\lambda \leq C (\lambda_{0,j} +\mu_0)(C_\eta
-{D_\eta}^{-1})N_j(\mu_0)\Big(1-(\mu_0 - \lambda_{0,j})^{-1} \Big). \] This is a
uniform estimate in $t$, $t\geq t_1$, and the right hand side goes to zero as
$\eta\rightarrow 0$. Combining these estimates, we conclude that
$I_{1}(t,\eta)\rightarrow 0$ {\em uniformly} in $t$, $t\geq t_1$ as $
\eta\rightarrow 0$. This completes the proof of the lemma.

\begin{lemma} : As $\eta$ goes to zero, $e^{\lambda_{0,j}t}\theta_{j,\nu}(t)$
converges uniformly to $e^{\lambda_{0,j}t}\theta_{j}(t)$ on $[0,+\infty[$.
\end{lemma}

\noindent{\bf Proof}: We estimate \[
\big|e^{\lambda_{0,j}t}\theta_{j,\nu}(t)-e^{\lambda_{0,j}t}\theta_{j}(t)\big| =
t\int^{\lambda_{0,j}+\nu}_{\lambda_{0,j}}e^{-(\lambda-\lambda_{0,j})t}N_{j}
(\lambda)d\lambda. \] Since the function $t\rightarrow
te^{-(\lambda-\lambda_{0,j})t}$ has its maximum at
$t=\frac{1}{\lambda-\lambda_{0,j}}\;$, we see that \[ \begin{array}{lcl}
\big|e^{\lambda_{0,j}t}\theta_{j,\nu}(t)-e^{\lambda_{0,j}t}\theta_{j}(t)\big|
\displaystyle & \leq &
\displaystyle\;e^{-1}\int^{\lambda_{0,j}+\nu}_{\lambda_{0,j}}\;\frac{ N_{j}
(\lambda)}{(\lambda-\lambda_{0,j})}\;d\lambda \\[+16pt] \displaystyle & \leq &
\displaystyle C\frac{e^{-1}}{\beta}{\nu_{\eta}} ^{\beta} \;\;\;\; \mbox{where
$\beta>0$, since ${\beta_{j}(M,g)}>0$ and by Lemma 1.2} . \end{array} \] This
proves the lemma.

\begin{lemma} : As $\eta$ goes to zero,
$e^{\lambda^{\widetilde{K},\widetilde{W}}_{0,j}t}\theta^{\widetilde{K},
\widetilde{W}}_{j,\nu}(t)-e^{\lambda^{\widetilde{K},\widetilde{W}}_{0,j}t}
\theta^{\widetilde{K},\widetilde{W}}_{j}(t)$ converges uniformly to zero on
$[0,+\infty[$. \end{lemma}

\noindent{\bf Proof}: We estimate

\[ \big|e^{\lambda^{\widetilde{K},\widetilde{W}}_{0,j}t}\theta^{\widetilde{K},
\widetilde{W}}_{j,\nu}(t)-e^{\lambda^{\widetilde{K},\widetilde{W}}_{0,j}t}
\theta^{\widetilde{K},\widetilde{W}}_{j}(t)\big|=t\int^{\lambda_{0,j}+\nu}_
{\lambda^{\widetilde{K},\widetilde{W}}_{0,j}}e^{-(\lambda-\lambda^{
\widetilde{K}, \widetilde{W}}_{0,j})t}N^{\widetilde{K},\widetilde{W}
}_{j}(\lambda)d\lambda. \]

By Proposition 2.10, $N^{\widetilde{K},\widetilde{W} }_{j}(\lambda)\leq
N_j(D_\eta\lambda)$, and because $\beta_j(M,g)>0$, \[N_j(\lambda)\leq
C(\lambda-\lambda_{0,j})^{\beta},\] for $\beta>0$ and $\lambda\leq
\lambda_{0,j}+1$. Therefore \[ \begin{array}{lcl}
\big|e^{\lambda^{\widetilde{K},\widetilde{W}}_{0,j}t}\theta^{\widetilde{K},
\widetilde{W}}_{j,\nu}(t)-e^{\lambda^{\widetilde{K},\widetilde{W}}_{0,j}t}
\theta^{\widetilde{K},\widetilde{W}}_{j}(t)\big| \displaystyle & \leq &
\displaystyle \;Cte^{\lambda^{\widetilde{K},
\widetilde{W}}_{0,j}t}\int^{\lambda_{0,j}+\nu}_{\lambda^{\widetilde{K},
\widetilde{W}}_{0,j}}\;e^{-\lambda
t}(D_\eta\lambda-\lambda_{0,j})^{\beta}\;d\lambda\\[+16pt] \displaystyle & \leq
& \displaystyle Cte^{\lambda^{\widetilde{K},
\widetilde{W}}_{0,j}t}(D_\eta(\lambda_{0,j}+\nu)-\lambda_{0,j})^{\beta}
\int^{\lambda_{0,j}+\nu}_{\lambda^{\widetilde{K},
\widetilde{W}}_{0,j}}\;e^{-\lambda t}\;d\lambda\\[+16pt] \displaystyle & \leq &
C(D_\eta(\lambda_{0,j}+\nu)-\lambda_{0,j})^{\beta} (1-e^{-\nu t}). \end{array}
\] This proves the lemma.

\noindent This also completes the proof of Theorem 4.1.

\bigskip

\noindent{\bf Remarks}. Suppose one could prove the following result,

\bigskip

{\it $\;$ Assume that $\beta_j(M, g) >0$. Given small $\varepsilon>0$, there are
positive constants $C_1$ and $C_2$ which are independent of the triangulation
such that} $$ C_1\,t^{-\bar\beta_j(M,g) - \varepsilon} \leq \theta_j(t) \leq
C_2\,t^{-\beta_j(M,g) + \varepsilon}. $$ \bigskip

\noindent Then, combined with Theorem 4.1, one could deduce a conjecture stated
in the introduction, that is,

\bigskip

{\it $\;$ Assume that $\beta_j(M, g) >0$. Then $\beta_j(K,g)$ converges to
$\beta_j(M,g)$ as the mesh of the triangulation goes to zero. Here
$\beta_j(K,g)$ denotes the combinatorial counterpart of $\beta_j(M,g)$.}

\noindent We are sadly as yet unable to improve the estimates in this section to
prove these conjectures.

\bigskip

\section{Calculations}

In this section, we calculate our von Neumann spectral invariants $
\beta_j(X,g)$ on closed hyperbolic manifolds, showing that they differ in
general from the Novikov-Shubin invariants. We also define Riemannian manifolds
$(M,g)$ with positive $\beta$-decay, that is, $\beta_j(M,g) >0$ is positive for
all $j$. We prove a result which gives some evidence to the conjecture stated in
the introduction.

Let $X$ be a closed hyperbolic manifold of dimension $d=2n+1$. This means that
$X =\Gamma\backslash G/K$ is a rank one locally symmetric space with $
G=SO_{0}(1,d)$, $K=SO (d)$ and $\Gamma\subset G$ is a torsion-free co-compact
discrete subgroup. We shall use some results of Fried \cite{F} to help us
compute our invariants for $X$. We first describe the Laplacian $ \Delta_{j}$ on
$j$-forms on $G/K$ in group theoretic terms as in \cite{F}.

It turns out that if one normalises the Killing form $c$ on the the Lie algebra
${\cal G}$ of $G$ to $\frac{1}{2d-2}c$, then the induced $G$ invariant metric on
$G/K$ has constant sectional curvature equal to $-1$.

Also the Casimir operator on $G$ induces a Casimir operator $\Omega_{j}$ on $
G/K$ acting on the space of $L^{2}$ sections of the homogeneous vector bundle of
$j$-forms on $G/K$. Then the Laplacian $\Delta_{j}=-\frac{1}{2d-2} \Omega_{j}$
is a constant multiple of the Casimir operator.

Consider the Iwasawa decomposition $G=KAN$ where $A$ has dimension one and hence
$G/K$ is a rank one symmetric space such that the rank of $G$ is greater than
the rank of $K$. Let $M$ be the centraliser of $A$ in $K$.

$G$ has no discrete series representations and the principal series
representations of $G$ are parametrised by $\hat M\times I\!\!R$, which carries
a smooth Plancherel density \cite{K}. On each line $\sigma\times I\!\!R$, it is
of the form $P_{\sigma}(\nu)d\nu$, where $P_{\sigma}(\nu)$ is an even polynomial
of degree $d-1$.

The $\sigma^{\prime}s$ which are of interest to us are those which occur in the
restriction of $\xi_{j}$ to the subgroup $M$, where $\xi_{j}$ denotes the usual
representation of $K=SO(d)$ on $\Lambda^{j}I\!\!\!\!C^{d}$. Writing
$I\!\!\!\!C^{d}=I\!\!\!\!C^{d-1}\oplus I\!\!\!\!C$, we observe that each
$\omega\in\Lambda^{j}I\!\!\!\!C^{d}$ is of the form $\omega^{\prime}+
\omega"\wedge dx_{d}$ where $\omega^{\prime}\in\Lambda^{j}I\!\!\!\!C^{d-1}$ and
$\omega"\in\Lambda^{j-1}I\!\!\!\!C^{d-1}$. Hence $\xi_{j}$ restricted to M is
isomorphic to $\sigma_{j}\oplus\sigma_{j-1}$, where $\sigma_{j}$ is the usual
representation of $M=SO(d-1)$ on $\Lambda^{j}I\!\!\!\!C^{d-1}$. Each $
\sigma_{j}$ is unitary and irreducible except in the case when $j=n$, in which
case it decomposes as a direct sum of two irreducible representations $
\sigma_{j}^{+}$ and $\sigma_{j}^{-}$.

The following theorem can be deduced from \cite{F}, theorem 2, observing that
the von Neumann trace of the heat kernel on $j$-forms on $G/K$ is just the
identity term in Fried's version of the Selberg Trace Formula for the trace of
the heat kernel on $j$-forms on $X$.

\begin{theorem} For $j=0,1,..,n$ we have \[ \theta_{j}(t) =
\tau(\exp(-t\Delta_{j})) = I_{t}(\sigma_{j}) + I_{t}(\sigma_{j-1}) \] where
$I_{t}(\sigma_{-1})\equiv 0$ and \[ I_{t}(\sigma_{j}) =
a_{j}\int^{\infty}_{-\infty}\exp(-t(\nu^{2} + c_{j}^{2})) P_{\sigma_{j}}(\nu)
d\nu \] Here $a_{j} = \left({d-1\atop j}\right) vol(X)$ and $c_{j}=n-j$.
\end{theorem}

By the isometry induced by the Hodge star operator, we see that \[ \theta_{j}(t)
= \theta_{d-j}(t) \] for $j=0,1,2,....,n$ and hence we obtain expressions for
$\theta_{j}(t)$ for $j=0,1,\ldots,d$. Using these, and the explicit expression
for the Plancherel measure \cite{K}, \cite{Mi}, we will be able to compute $
\beta_j(X)$.

\begin{theorem} . Let $X$ be a closed hyperbolic manifold of dimension $d=2n+1$.
Then $ {\overline\beta_j(X)} = \beta_j(X) = \frac{3}{2} $ for $j=0,1,..,n-1$ and
$ {\overline\beta_n(X)} = \beta_n(X) = \frac{1}{2} $ . \end{theorem}

\noindent{\bf Proof}. Using the following explicit expression for the Plancherel
measure, \[ P_{\sigma_{j}}(\nu) = \Big(\nu^2 +
(n-j)^2\Big)^{-1}\prod_{k=0}^n\Big(\nu^2 + k^2\Big) \] we see that for $j\leq
n-1$, \[ \begin{array}{l} \displaystyle I_{t}(\sigma_{j}) =
C_1\exp(-tc_{j}^{2})\Big\{ \int^{\infty}_{-\infty} \exp(-t\nu^{2})\nu^2 d\nu +
O(t^{-{\frac{5}{2}}}) \Big\} \\ \\ \displaystyle =
C_1^{\prime}\exp(-tc_{j}^{2})\Big\{t^{-{\frac{3}{2}}} +
O(t^{-{\frac{5}{2}}})\Big\}\;\;\;\; {\rm as}\;\;\; t\rightarrow\infty.
\end{array} \] Using the previous theorem, we see that ${\overline\beta_j(X)} =
\beta_j(X) = \frac{3}{2} $ for $ j=0,1,..,n-1$.

The case when $j=n$ remains to be studied. In this case \[ I_{t}(\sigma_{n}) =
C_2\int^{\infty}_{-\infty} \exp(-t\nu^{2}) d\nu + O(t^{-{
\frac{5}{2}}})\\\\=C_2^{\prime}t^{-{\frac{1}{2}}}+ O(t^{-{\frac{3}{2}}}). \]
Using the previous theorem, we see that ${\overline\beta_n(X)} = \beta_n(X) =
\frac{1}{2} $.

\begin{definition} . A closed Riemannian manifold $(M,g)$ is said to have
positive $\beta$-decay if $\beta_j(M,g) >0$ is positive for all $j$.
\end{definition}

We recall that a closed Riemannian manifold $(M,g)$ is said to be $L^2$ -acyclic
if all the $L^2$ Betti numbers of its universal cover vanish, that is,
$b^j_{(2)}(\widetilde M) = 0$ for all $j$. (cf. \cite{M1}, \cite{CM}.)

\begin{proposition} . Let $(M,g)$ be a closed $L^2$ acyclic Riemannian manifold
with positive $ \beta$-decay, and $(N,h)$ be any closed $L^2$ acyclic Riemannian
manifold. Then $(M\times N, g\times h)$ has positive $\beta$-decay.
\end{proposition}

\noindent{\bf Proof}. Since $M\times N$ is given the product metric $g\times h$,
\[ \tau_{M\times N}\Big(e^{-t\Delta_k^{M\times N}}\Big) = \sum_{i+j=k} \tau_{M}
\Big(e^{-t\Delta_i^{M}}\Big) \tau_{N}\Big(e^{-t\Delta_j^{N}}\Big) \] and \[
\lambda_{0,k}^{M\times N} = \mbox{minimum}\Big\{\lambda_{0,i}^{M}+
\lambda_{0,j}^{N} : i+j =k\Big\}. \] Since $M$ and $N$ are $L^2$-acyclic \[
\theta_k^{M\times N}(t) = \sum_{i+j=k}\theta_i^{M}(t)\theta_j^{N}(t). \] Hence
we see that \[ \beta_k(M\times N, g\times h) = \mbox{minimum}\Big\{\beta_i(M, g)
+ \beta_j(N, h) : i+j =k\Big\} \] The theorem follows.

\begin{corollary} . Let $(M,g)$ be a closed, odd dimensional, hyperbolic
manifold with positive $\beta$-decay, and $(N,h)$ be any closed $L^2$ acyclic
Riemannian manifold. Then $(M\times N, g\times h)$ has positive $\beta$-decay.
\end{corollary}

\bigskip

\section{Von Neumann Determinants and $\beta$-Torsion}

In this section, we assume that $(M,g)$ is a closed $n$-dimensional Riemannian
manifold such that $\beta_j(M,g)>0$. We will define the von Neumann determinant
of the operator $\Delta_j - \lambda_{0,j}$ on $ \widetilde{M}$ following ideas
of \cite{M1} and \cite{L}, where the von Neumann determinant of the operator
$\Delta_j$ was studied. Hence we will also assume that $\lambda_{0,j} > 0$,
without any loss of generality. We will compute these determinants for certain
closed hyperbolic dimensional manifolds. We also define the analytic
$\beta$-torsion in terms of the von Neumann determinants of the operators
$\Delta_j - \lambda_{0,j}$, by analogy to Ray-Singer torsion \cite{RS}. We also
define its combinatorial counterpart, which we call combinatorial
$\beta$-torsion. In the next section, we prove a result which gives evidence
that the combinatorial $\beta$-torsion converges to the analytic
$\beta$-torsion, as the mesh of the triangulation goes to zero.

We begin by defining the partial $L^{2}$ zeta functions of the operator
$\Delta_j - \lambda_{0,j}$ as follows.

\begin{definition} . The partial $L^{2}$ zeta functions of the operator
$\Delta_j - \lambda_{0,j}$ are \[ \begin{array}{rcl} \zeta^{(1)}_{j}(s) & = &
\displaystyle\frac{1}{\Gamma(s)} \int^{1}_{0}t^{s-1} e^{\lambda_{0,j}
t}\theta_{j}(t)dt \\[+16pt] \mbox{and }\;\;\zeta^{(\infty)}_{j}(s) & = &
\displaystyle\frac{1}{ \Gamma(s)} \int^{\infty}_{1}t^{s-1}e^{\lambda_{0,j} t}
\theta_{j}(t) dt. \end{array} \] Here $s$ belongs to a subset of complex numbers
which will be specified shortly. \end{definition}

We begin by proving the following lemma.

\begin{lemma} . $\zeta^{(1)}_{j}(s)$ is a holomorphic function in the half-plane
$ \Re(s)>n/2$ and has a meromorphic continuation to $I\!\!\!C$ with no pole at $
s=0$. \end{lemma}

\noindent{\bf Proof}. Using \cite{R}, we have an asymptotic expansion as $
t\rightarrow 0^{+}$ of $\tau(e^{-t\Delta_j})$, \[ \tau(e^{-t\Delta_j})\sim
t^{-n/2}\sum_{j=0}^{\infty}t^{j}c_{j}\;\;\mbox{ as } \;\; t\rightarrow 0^{+}. \]
Hence $e^{\lambda_{0,j} t}\tau(e^{-t\Delta_j})$ has the following small time
asymptotic expansion \[ e^{\lambda_{0,j} t}\tau(e^{-t\Delta_j})\sim
t^{-n/2}\sum_{i=0}^{\infty}\sum_{l+k=i}c_{l}{\frac{\lambda_{0,j}^k}{k !}}
t^{i}\;\;\mbox{ as }\;\; t\rightarrow 0^{+}. \] In particular, $e^{\lambda_{0,j}
t}\tau(e^{-t\Delta_j})\leq Ct^{-n/2}$ for $ 0\leq t\leq 1$. We deduce that
$\zeta^{(1)}_{j}(s)$ is well defined on the half-plane $\Re(s)>n/2$. Clearly \[
\frac{\partial}{\partial\bar{s}}\;\zeta^{(1)}_{j}(s)=0 \] in this region, that
is, $\zeta^{(1)}_{j}(s)$ is holomorphic in this half-plane.

The meromorphic continuation of $\zeta^{(1)}_{j}(s)$ to the half-plane $
\Re(s)>n/2-N$ is obtained by considering the first $N$ terms of the small time
asymptotic expansion of $e^{\lambda_{0,j} t}\tau(e^{-t\Delta_j})$, \[
\begin{array}{lcl} \zeta^{(1)}_{j}(s) & = &
\displaystyle-\;\Big(\frac{e^{\lambda_{0,j}}-1}{ \lambda_{0,j}}\Big)
\frac{b^{j}_{(2)} (\widetilde M)}{s\Gamma(s)}+ \frac{1}{ \Gamma(s)}
\left[\int^{1}_{0}t^{s-1-n/2}\Big(\sum_{i=0}^{N}\sum_{l+k=i}c_{l}{
\frac{\lambda_{0,j}^k}{k !}}t^{i}\Big)dt+R_{N}(s)\right] \\[+16pt] & = &
\displaystyle -\;\Big(\frac{e^{\lambda_{0,j}}-1}{\lambda_{0,j}}\Big)
\frac{b^{j}_{(2)}(\widetilde M)}{s\Gamma(s)}+\frac{1}{\Gamma(s)}\left[
\sum_{i=0}^{N}\;\frac{\sum_{l+k=i}c_{l}{\frac{\lambda_{0,j}^k}{k !}}}{
(s+i-n/2)}+R_{N}(s)\right] \end{array} \] where $R_{N}(s)$ is holomorphic in the
half plane $\Re(s)>n/2-N$.

By observation, we see that the meromorphic continuation of $\zeta^{(1)}_{j}
(s)$, also denoted by $\zeta^{(1)}_{j}(s)$, has no pole at $s=0$.

\begin{corollary} . $\zeta^{(1)}_{j}(0)=\left\{ \begin{array}{lcl}
-\Big(\frac{e^{\lambda_{0,j}}-1}{\lambda_{0,j}}\Big) b^{j}_{(2)}(\widetilde M)+
\sum_{l+k={n\over 2}}c_{l}{\frac{\lambda_{0,j}^k}{k !}}\;\; & \mbox{ if } & n
\mbox{ is even} \\ -\Big(\frac{e^{\lambda_{0,j}}-1}{\lambda_{0,j}}\Big)
b^{j}_{(2)}(\widetilde M) & \mbox{ if } & n\mbox{ is odd}. \end{array} \right.$
\end{corollary}

\begin{lemma} . $\zeta^{(\infty)}_{j}(s)$ is holomorphic in the half-plane $
\Re(s)<\beta_j(M,g)$. \end{lemma}

\noindent{\bf Proof}. Using the estimate $e^{\lambda_{0,j}t}\theta_{j}(t) \leq
Ct^{-\beta_j(M,g)+\varepsilon}$ we see that $\zeta^{(\infty)}_{j}(s)$ is well
defined on the half-plane $\Re(s)<\beta_j(M,g)$, since \[
\displaystyle\Big|\zeta^{(\infty)}_{j}(s)\Big|  \leq \displaystyle\frac{C}{
|\Gamma(s)(\Re(s)-\beta_j(M,g))|} \] whenever $\Re(s)<\beta_j(M,g)$. Clearly
$\displaystyle\frac{\partial}{ \partial \bar s}\;\zeta^{(\infty)}_{j}(s)=0$ on
this half-plane.

The following is an immediate consequence of the proof of Lemma 6.4.

\begin{corollary} . $\zeta^{(\infty)}_{j}(0)=0$. \end{corollary}

\begin{definition} . Define the $L^{2}$ zeta function of the operator $\Delta_j
- \lambda_{0,j}$ as follows. \[
\zeta_{j}(s)=\zeta^{(1)}_{j}(s)+\zeta^{(\infty)}_{j}(s). \] \end{definition}

The following theorem summarizes the prior lemmas.

\begin{theorem} . $\zeta_{j}(s)$ is holomorphic near $s=0$. Its value at zero is
\[ \zeta_{j}(0) = \left\{ \begin{array}{lcl}
-\Big(\frac{e^{\lambda_{0,j}}-1}{\lambda_{0,j}}\Big) b^{j}_{(2)}(\widetilde
M)+\sum_{l+k={n\over 2}} c_{l}{\frac{\lambda_{0,j}^k}{k !}}\;\; & \mbox{ if } &
n \mbox{ is even} \\ -\Big(\frac{e^{\lambda_{0,j}}-1}{\lambda_{0,j}}\Big)
b^{j}_{(2)}(\widetilde M) & \mbox{ if } & n\mbox{ is odd}. \end{array} \right.
\] \end{theorem}

\begin{definition} . The von Neumann determinant of the operator $\Delta_j -
\lambda_{0,j}$ on $ \widetilde{M}$ is by definition \[
|\mbox{Det}_{\tau}|(\Delta_j - \lambda_{0,j})=\exp(-\zeta^{\prime}_{j}(0)). \]
\end{definition}

We now compute this determinant on hyperbolic manifolds. In principle, it is
possible to use a computer program and compute all the determinants. As an
illustration, we compute the von Neumann determinant of $\Delta_1 -
\lambda_{0,1} = \Delta_1 - 1$ on five dimensional hyperbolic space. We use the
notation of section 5.

Since there are no $L^2$ harmonic differential forms on five dimensional
hyperbolic space, we see that $\zeta_j(0) = 0$ by Theorem 6.7. By Theorem 5.1,
\[ \theta_{1}(t) = \tau(\exp(-t\Delta_{1})) = I_{t}(\sigma_{1}) +
I_{t}(\sigma_{0}) \] where \[ I_{t}(\sigma_{j}) =
a_{j}\int^{\infty}_{-\infty}\exp(-t(\nu^{2} + c_{j}^{2})) P_{\sigma_{j}}(\nu)
d\nu \] Here $j=0,1$, $a_{j} = \left({4\atop j}\right) vol(X)$ and $
\lambda_{0,j} = c_{j}^2= (2-j)^2$.

An easy calculation yields

\begin{lemma} . The von Neumann determinant of $\Delta_1 - \lambda_{0,1}$ is
given by \[ \log |\mbox{Det}_{\tau}|(\Delta_1 - 1) = 8.7062\;{\mbox{vol}}(X) \]
\end{lemma}

Next we define the analytic $\beta$-torsion to be

\begin{definition}. Suppose that $M$ is a manifold with positive $\beta$-decay.
Then the analytic $\beta$-torsion is defined to be the product $$ T(\widetilde
M, g) = \prod_{j=0}^n |\mbox{Det}_{\tau}|(\Delta_j - \lambda_{0,j})^{j(-1)^j} $$
\end{definition}

One can make an analogous definition in the combinatorial setting. The results
of the next section suggest that the combinatorial torsion may converge to its
analytic counterpart. We plan to discuss this further elsewhere.

\bigskip

\section{Convergence of Theta and Zeta Functions}

Using the results of section 3, we prove in this section that the $L^2$
combinatorial theta function, as a function of the time variable converges
uniformly on compact subsets to the $L^2$ analytic theta function as the mesh of
the triangulation goes to zero.

Then we are able to prove the analogue of this result for $L^2$-zeta functions.
However, to define these zeta functions, one needs to impose the condition of
positive decay on the manifold. Our method of circumventing the slow decay of
the $L^2$ theta functions is to split the $L^2$ zeta functions into two partial
zeta functions. We introduce an analogous splitting of the combinatorial $L^2$
zeta function and prove that the resulting partial zeta functions converge as
the mesh of the triangulation goes to zero to the corresponding $L^2$ analytic
partial zeta functions, uniformly on compact subsets of their respective
domains. These results are analogs for our situation of the Dodziuk-Patodi
theorems \cite{DP}.

If our manifold has positive $\beta$-decay, then using Theorem 4.1, we prove
that the zeta functions of the operator $\Delta_j ^{\widetilde{K},\widetilde{W}}
-\lambda_{0,j}^{\widetilde{K},\widetilde{W}}$ converge to the zeta functions of
the operator $\Delta_j -\lambda_{0,j}$ as the mesh of the triangulation goes to
zero. Recall that we have also come across these zeta functions in the previous
section.

We begin with

\begin{theorem} : As the mesh goes to zero,
$\;\theta_{j}^{\widetilde{K},\widetilde{W}}(t)$ converges to $\theta_{j}(t)$
uniformly on the set $[t_1,\infty)$ for any $ t_1>0$. \end{theorem}

\noindent{\bf Proof}: $\;$By Propositions 3.1 and 3.2, we see that \[
\left|\theta^{\widetilde{K},\widetilde{W}}_j(t)-\theta_{j}
(t)\right|\leq\left|\theta_{j}(D^{-1}_{\eta}t)-\theta_{j}(t) \right|+
\left|\theta_{j}(C_\eta t)-\theta_{j}(t)\right|+ \varepsilon(\eta). \] So it
suffices to prove that $\theta_{j}(B_\eta t)$ converges, uniformly in $
t\in[t_1,\infty)$, to $\theta_{j}(t)$ as $\eta\rightarrow 0$. Here $
B_{\eta}\rightarrow 1^+$ as $\eta\rightarrow 0$. By the mean value theorem, we
see that for $t\geq t_1$, \[
|\theta_{j}(B_{\eta}t)-\theta_{j}(t)|\leq|B_{\eta}-1|M \] where
$\displaystyle\;M=\sup_{t\geq t_1}|t\;\tau(\Delta_{j}e^{-t\Delta_{j}})|$. It
suffices to show that $ M$ is finite.

Now, \[ \begin{array}{lcl} \displaystyle t\;\tau(\Delta_{j}e^{-t\Delta j}) & = &
\displaystyle t\int^{\infty}_{0}\lambda e^{-t \lambda}dN_{j}(\lambda) \\[+14pt]
\displaystyle & = & \displaystyle
-t\int^{\infty}_{0}e^{-t\lambda}N_{j}(\lambda)d\lambda+t^{2}
\int^{\infty}_{0}\lambda e^{-t\lambda}N_{j}(\lambda)d\lambda \\[+14pt]
\displaystyle & = & \displaystyle
-\tau(e^{-t\Delta_{j}})+\int^{\infty}_{0}\lambda e^{-\lambda}
N_{j}\left(\frac{\lambda}{t}\right)d\lambda. \end{array} \] So for $t\geq t_1$,
\[ |t\;\tau(\Delta_{j}e^{-t\Delta_{j}})|\leq b^{j}_{(2)}(\widetilde{M}
)+\theta_{j}(t_1)+\int^{\infty}_{0}\lambda e^{-\lambda}N_{j}\left(\frac{
\lambda} {t_1}\right)d\lambda. \] That is, $M<\infty$.

Let $M$ denote a closed manifold. Recall the definition of the Novikov-Shubin
invariants, \[ \alpha_{j}(M)=\sup\{\beta\in I\!\!R: \theta_{j}(t)\;\;\mbox{ is }
\;\;O(t^{-\beta})\;\mbox{ as }\;t\rightarrow\infty\}\in[0,\infty]. \] It has
been shown by Novikov and Shubin \cite{ES} that the numbers $ \alpha_{j}(M)$ are
independent of the choice of metric on $M$, and in fact Gromov and Shubin
\cite{GS} show that they depend only on the homotopy type of $M$. We recall from
\cite{M1} the following definition.

\begin{definition} : A manifold $M$ is said to have positive decay if all of its
Novikov-Shubin invariants $\alpha_{j}(M)$ are positive, for $j\geq 0$.
\end{definition}

Let $K$ be a smooth triangulation of $M$. Recall the definition of the
combinatorial analogue of the Novikov-Shubin invariants, \[
\alpha_{j}(K)=\sup\{\beta\in I\!\!R: \theta^{\widetilde{K},\widetilde
W}_{j}(t)\;\;\mbox{ is } \;\;O(t^{-\beta})\;\mbox{ as }\;t\rightarrow\infty\}
\in[0,\infty]. \] It has been shown by A.V. Efremov \cite{E} that the numbers
$\alpha_{j}(K)$ are independent of the choice of triangulation of $M$, and hence
they depend only on the topology of $M$. In fact, Efremov proves that \[
\alpha_{j}(K) = \alpha_{j}(M). \] Note that one can extract the proof of his
result from sections 2, 3 and 4.

See \cite{M1,M2}, \cite{L} for numerous examples of manifolds with positive
decay.

In the following discussion, we will only consider manifolds with positive
decay.

Now we define the zeta function of the Laplacian $\Delta_j$ following \cite
{M1}, and recall some of its properties. Introduce two partial zeta functions:
$$ \zeta_{j}(s,1)= {\frac{1}{\Gamma(s)}}\int_0^{1} \theta_{j}(t)t^{s-1}dt, \ \ \
\ (\Re(s)> {\frac{1}{2}}{\dim M}), \eqno $$

and under the assumption that $M$ has positive decay, $$ \zeta_{j}(s,\infty) =
{\frac{1}{\Gamma(s)}}\int_{1}^\infty \theta_{j}(t)t^{s-1}dt \ \ \ \
(\Re(s)<\alpha_j(M)). \eqno $$

Then $\zeta_{j}(s,\infty)$ is analytic for $\Re(s)<\alpha_j(M)$ whereas $
\zeta_{j}(s,1)$ is analytic for $\Re(s) > {\frac{1}{2}}{\dim M}$ and has an
analytic continuation to a neighbourhood of zero. Thus with this understanding
we can define the zeta function to be the sum $$
\zeta_{j}(s)=\zeta_{j}(s,\infty)+\zeta_{j}(s,1). \eqno $$

It is analytic on a neighbourhood of zero.

We have seen that the assumption of positive decay implies decay for $
\theta_{j}^{\widetilde K, \widetilde W}(t)$ (with the same $\alpha_j(M)$) and
hence we may define zeta functions $\zeta_{j}^{\widetilde K, \widetilde W}(s)$
similarly in terms of $\theta_{j}^{\widetilde K,\widetilde W}(t)$ (note however
that in the combinatorial case the theta functions are well defined at $t=0$ and
so the analytic continuation we employed above is not necessary).

\begin{theorem}. Let $M$ be a closed manifold with positive decay. Then \\ 1.
the combinatorial partial zeta function \[ \zeta_{j}^{\widetilde{K},\widetilde
W}(s,1)\equiv\frac{1}{\Gamma(s)} \;\int^{1}_{0}t^{s-1}
\theta_{j}^{\widetilde{K},\widetilde{W}}(t)dt \] converges for
$\Re(s)>{\frac{1}{2}}{\dim M} = {n\over 2}$ uniformly on compact subsets to \[
\zeta_{j}(s,1)\equiv\frac{1}{\Gamma (s)}\;\int^{1}_{0}t^{s-1}\theta_ {j}(t)dt \]
2. while the other combinatorial partial zeta function: \[
\zeta_{j}^{\widetilde{K},\widetilde W}(s,\infty)=\frac{1}{\Gamma(s)}
\;\int^{\infty}_{1}t^{s-1} \theta^{\widetilde{K},\widetilde W}_{j}(t)dt \]
converges for $\Re(s)<\alpha(M) = {\rm min}\{\alpha_j(M)\}$ uniformly on compact
subsets to \[
\zeta_{j}(s,\infty)\equiv\frac{1}{\Gamma(s)}\;\int^{\infty}_{1}t^{s
-1}\theta_{j}(t) dt. \] \end{theorem}

\noindent{\bf Proof}.  {\em Part 1}.  Let $F$ be a compact subset of the
half-plane $\Re s > {n\over 2}$, that is $\Re s \geq {n\over 2} + \delta$ for
all $s\in F$. We need to prove that given $\varepsilon > 0$ $$
\int^{1}_{0}dt\left|\int^{\infty}_{0}t^{s}e^{-t\lambda}\left(N_{j}(\lambda)-N^
{\widetilde{K},\widetilde{W}}_{j}(\lambda)\right)d\lambda\right|<\varepsilon. $$
for all $\eta$ small enough and for any $s\in F$.

Since $\theta_j(t) \leq C t^{-{n\over 2}}$ for $0<t\leq 1$, one sees that $$
|\Gamma(s)\zeta_{j}(s,1)|\leq\int^{1}_{0}dt\int^{\infty}_{0}t^{\Re(s)}e^{-t
\lambda}N_{j}(\lambda)d\lambda<\infty $$ is finite and bounded for any $s\in
F$.Therefore, given $\varepsilon>0$, there is an $m>0$ such that $$
\int^{1}_{0}dt\int^{\infty}_{m}t^{\Re(s)}e^{-t\lambda}N_{j}(\lambda)d\lambda
<\varepsilon. $$ By the uniform convergence theorem and using Propositions 2.2
and 2.6, one has for $\eta$ small enough, $$
\int^{1}_{0}dt\left|\int^{m}_{0}t^{s}e^{-t\lambda}\left(N_{j}(\lambda)-N^
{\widetilde{K},\widetilde{W}}_{j}(\lambda)\right)d\lambda\right|<\varepsilon $$
uniformly for $s\in F$. In the notation of Proposition 2.2, define the interval
$$ S_{m}=\Big(m,\Big\{{1\over C_1{\sqrt h}} - 1\Big\}^2\Big) $$ where $h =
(\eta|\log\eta|)^{2}$. Then by Proposition 2.2, for all $\lambda\in S_{m}$ one
has the inequality $$ N^{\widetilde{K},\widetilde{W}}_{j}(\lambda)\leq
N_{j}(D_\eta\lambda). $$ On the other hand, if $\lambda>m$ and $\lambda\not\in
S_{m}$, then there is a $\beta>0$ such that $$
\lambda\geq\beta(\eta|\log\eta|)^{-2} $$ for all $\eta$ small enough.  By the
fullness assumption on $K$, one has $$
N^{\widetilde{K},\widetilde{W}}_{j}(\lambda)\leq\gamma\eta^{-n} $$ for some
$\gamma>0$.

One estimates for $s\in F$,

$$ \begin{array}{lcl} \displaystyle
\int^{1}_{0}dt\left|\int^{\infty}_{m}t^{s}e^{-\lambda
t}N^{\widetilde{K},\widetilde{W}}_ {j}(\lambda)d\lambda\right| & \leq &
\displaystyle\int^{1}_{0}dt\left|\int_{\lambda\in S_{m}}t^{\Re(s)}e^{-\lambda t}
N_{j}(D_\eta\lambda)d\lambda\right| \\[+18pt]
&+&\displaystyle\;\int^{1}_{0}dt\left|\int_{\beta\over h}^\infty
t^{\Re(s)}e^{-\lambda
t}N^{\widetilde{K},\widetilde{W}}_{j}(\lambda)d\lambda\right|\\[+18pt] &\leq &
D^{\Re(s)}_\eta\varepsilon+\gamma\beta^{-\Re(s)}\eta^{2\delta}|\log\eta|^{2\
Re(s)}\;. \end{array} $$ Thus there is a constant $C>0$ such that for $s\in F$,
one has $$ \int^{1}_{0}dt\left|\int^{\infty}_{m}t^{s}e^{-\lambda
t}N^{\widetilde{K},\widetilde{W}}_ {j}(\lambda)d\lambda\right|\leq C\varepsilon
$$ for all $\eta$ small enough.

This proves Part 1.

\noindent{\em Part 2}.  Let $F$ be a compact subset of the half-plane $\Re(s) <
\alpha_j(M)$, that is $\Re(s) \leq \alpha_j(M) + \delta$ for all $s\in F$. We
need to prove that given $\varepsilon > 0$ $$ \int_1^\infty t^{{\Re(s)} - 1}
\left|\theta^{\widetilde{K},\widetilde W}_{j}(t) - \theta_{j}(t)\right| dt<
\varepsilon $$ for all $\eta$ small enough. For $s\in F$, we see as in the first
line of the proof of Theorem 7.1, that the sequence of functions $$ t^{{\Re(s)}
- 1} \left|\theta^{\widetilde{K},\widetilde W}_{j}(t) - \theta_{j}(t)\right| $$
is dominated, for $t\in [1,\infty)$, by the function $$ t^{{\Re(s)} -
1}\left|\theta_{j}({t\over 2})-\theta_{j}(t) \right|+ t^{{\Re(s)} -
1}\left|\theta_{j}(2 t)-\theta_{j}(t)\right|+ t^{{\Re(s)} -
1}\varepsilon(t).\eqno(*) $$ The first two terms above are clearly integrable on
the interval $[1,\infty)$. We now examine the last term. From section 3, we see
that $$ t^{{\Re(s)} - 1}\varepsilon(t) \leq 2 t^{{\Re(s)}}\int_\mu^\infty
e^{-2t\lambda}N_j(\lambda) d\lambda + {1\over 2} t^{{\Re(s)}}\int_\mu^\infty
e^{-{t\lambda \over 2}}N_j(\lambda) d\lambda. $$ We observe that for $a>0$, one
has $$ \int_1^\infty t^{{\Re(s)}}\int_\mu^\infty e^{- a t\lambda}N_j(\lambda)
d\lambda \leq \int_1^\infty t^{{\Re(s)}-1} \theta_j(a t) dt. $$ Therefore for
$s\in F$, one has $$ \int_1^\infty t^{{\Re(s)} - 1}\varepsilon(t) dt \leq 2
\int_1^\infty t^{{\Re(s)}-1} \theta_j(2 t) dt + {1\over 2}\int_1^\infty
t^{{\Re(s)}-1} \theta_j({t\over 2}) dt<\infty $$ and the third term in $(*)$ is
integrable on the interval $[1,\infty)$. By Theorem 7.1, we see that $$
\lim_{\eta\rightarrow 0}\,\,t^{{\Re(s)}-1}
\left|\theta^{\widetilde{K},\widetilde W}_{j}(t) - \theta_{j}(t)\right| = 0 $$
pointwise, for any $t\in [1,\infty)$. By the Dominated Convergence Theorem, we
conclude that for $s\in F$ $$ \int_1^\infty t^{{\Re(s)} - 1}
\left|\theta^{\widetilde{K},\widetilde W}_{j}(t) - \theta_{j}(t)\right| dt<
\varepsilon $$ for all $\eta$ small enough, proving part 2.

This is the analogue for $L^2$-zeta functions of the results of \cite{DP} for
the ordinary zeta function. Notice that in fact the proofs make no essential use
of the representation of the fundamental group (they use only the fact that the
commutant is a semi-finite von Neumann algebra).

For the next theorem, we assume that the following conjecture holds,

\bigskip

\noindent{\em $\;$ Assume that $\beta_j(M, g) >0$. Then $\beta_j(K,g)$ converges
to $\beta_j(M,g)$ as the mesh of the triangulation goes to zero. Here
$\beta_j(K,g)$ denotes the combinatorial counterpart of $\beta_j(M,g)$.}

\bigskip

This conjecture was stated in the introduction and also discussed in section 4.
Then in the notation of section 6, Theorem 4.1 combined with dominated
convergence implies:

\begin{theorem}. Let $(M,g)$ be a closed Riemannian manifold with positive
$\beta$-decay. Then the combinatorial partial zeta function of the operator
$\Delta_j^{\widetilde K, \widetilde W}-\lambda_{0,j}^{\widetilde K, \widetilde
W}$ \[ \zeta_{j}^{\widetilde{K},\widetilde W (1)}(s)\equiv\frac{1}{\Gamma(s)}
\;\int^{1}_{0} t^{s-1}
e^{\lambda^{\widetilde{K},\widetilde{W}}_{0,j}t}\theta_{j}^{\widetilde{K},
\widetilde{W}}(t)dt
\] converges for $\Re(s)>{\frac{1}{2}}{\dim M}$ uniformly on compact subsets to
\[ \zeta_{j}^{(1)}(s)\equiv\frac{1}{\Gamma (s)}\;\int^{1}_{0}
t^{s-1}e^{\lambda_{0,j}t}\theta_ {j}(t)dt \] while the other combinatorial
partial zeta function of the operator $\Delta_j^{\widetilde K, \widetilde
W}-\lambda_{0,j}^{\widetilde K, \widetilde W}$: \[
\zeta_{j}^{\widetilde{K},\widetilde W (\infty)}(s)=\frac{1}{\Gamma(s)}
\;\int^{\infty}_{1} t^{s-1} e^{\lambda^{\widetilde{K},\widetilde{W}}_{0,j}t}
\theta^{\widetilde{K},\widetilde W}_{j}(t)dt \] converges for $\Re(s)<\beta=
{\rm min}\{\beta_j(M,g)\}$ uniformly on compact subsets to \[
\zeta_{j}^{(\infty)}(s)\equiv\frac{1}{\Gamma(s)}\;\int^{\infty}_{1}
t^{s-1}e^{\lambda_{0,j}t}\theta_{j}(t) dt. \] \end{theorem}

\noindent{\bf Proof}.The proof is similar to that given in the previous theorem,
but now using the results in sections 2 and 4 instead. The only new point to
observe is that the assumption that the conjecture holds is used to show that
when the mesh of the triangulation is small enough, one sees that $\beta_j(K,g)$
is positive and so the combinatorial partial zeta function
$\zeta_{j}^{\widetilde{K}, \widetilde W (\infty)}(s)$ is defined. The detailed
proof will be omitted.

\bigskip

\section{Appendix}

Here we present a proof due to Terry Lyons that $\bar\beta_{0}(M,g) \geq 1$ and
$\beta_{0}(M,g) \geq 1$. He has kindly permitted us to do so.

\noindent{\bf Proposition}. {\em $\bar\beta_{0}(M,g) \geq 1$ and $\beta_{0}(M,g)
\geq 1$.}

\noindent{\bf Proof}.  Let $\lambda_{0,0}$ be the spectral gap.  Then there is a
$v>0$, $v\in L^{\infty}(\widetilde{M})$ such that $$ \Delta_{0}v=\lambda_{0,0}
v. $$ That is, $v$ is a ground state for $\Delta_{0}$ which is unique up to a
multiplicative constant.  Then either there is a $\gamma\in\pi_{1}(M)$ such that
$\gamma^\star v$ is not proportional to $v$, or for every $\gamma\in\pi_{1}(M)$,
$$ \gamma^\star v=\varphi(\gamma)v $$ for some morphism
$\varphi:\pi_{1}(M)\rightarrow I\!\!R_{+}$.

\noindent{\bf Case 1}.  Suppose that there is a $\gamma\in\pi_{1}(M)$ such that
$\gamma^\star v$ is not proportional to $v$.  Then $$ u=\frac{\gamma^\star
v}{v}>0 $$ is a positive, non-constant harmonic function for the Markovian
semi-group with integral kernel $$ \tilde{p}_{t}(x,y)=e^{\lambda_{0,0}
t}p_{t}(x,y)\;\frac{v(y)}{v(x)}\;. $$ Since 1 and $u$ are non-proportional
positive harmonic functions for this semigroup, it is a classical result then
that $\tilde{p}_{t}$ is transient cf.\cite{An} page 44, that is, $$
\int^{\infty}_{1}e^{\lambda_{0,0} t}p_{t}(x,x)dt<\infty. $$

\noindent{\bf Case 2}.  Suppose that for every $\gamma\in\pi_{1}(M)$, $$
\gamma^\star v=\varphi(\gamma)v $$ for some morphism
$\varphi:\pi_{1}(M)\rightarrow I\!\!R_{+}$.  Then the Markovian semigroup $$
\tilde{p}_{t}(x,y)=e^{\lambda_{0,0} t}p_{t}(x,y)\;\frac{v(y)}{v(x)} $$ is
$\pi_{1}(M)$ invariant, by a simple calculation.  Now by \cite{LS} theorem 3,
$\tilde{p}_{t}$ admits a non-constant, positive harmonic function $u$. So we see
again that $\tilde{p}_{t}(x,y)$ is transient.

We conclude that $\tilde{p}_{t}(x,y)$ is always transient. It follows
immediately that $\bar\beta_{0}(M,g) \geq 1$. To prove that $\beta_{0}(M,g) \geq
1$, we further argue as follows.

Let $K = \displaystyle \int^{\infty}_{1}e^{\lambda_{0,0} t}p_{t}(x,x)dt<\infty.$
Recall that the function $t\rightarrow e^{\lambda_{0,0} t}p_{t}(x,x)$ is
non-increasing. Therefore $$ {t\over 2}\,e^{\lambda_{0,0} t}p_{t}(x,x)\leq
\int^{t}_{t \over 2}e^{\lambda_{0,0} s}p_{s}(x,x)ds \leq K $$ and
$\beta_{0}(M,g) \geq 1$.

\end{document}